\documentclass[smallcondensed]{svjour3}
\usepackage[T1]{fontenc}
\usepackage[latin9]{inputenc}
\usepackage{url}
\usepackage{amsmath}
\usepackage{stackrel}
\usepackage{graphicx}
\usepackage{epstopdf}
\makeatletter
\RequirePackage{fix-cm}
\usepackage[numbers]{natbib}
\smartqed  

\makeatother

\usepackage{babel}
\begin{document}
\title{Numerical solution for the stress near a hole with corners in an infinite
plate under biaxial loading
\thanks{This work was supported by a grant from the Simons Foundation (Award \#354717, BJS).}}
\author{Weiqi Wang \and Brian J. Spencer}
\institute{Weiqi Wang \at Department of Mathematics, University at Buffalo,
Buffalo, NY 14260, USA\\
Tel.: +(716)-352-3644\\
\email{weiqiwan@buffalo.edu}\\
\and Brian J. Spencer \at Department of Mathematics, University
at Buffalo, Buffalo, NY 14260, USA\\
Tel.: +(716)-645-8805\\
\email{spencerb@buffalo.edu}\\
}
\date{Received: date / Accepted: date}
\maketitle
\begin{abstract}
We consider the elastic stress near a hole with corners in an infinite
plate under biaxial stress. The elasticity problem is formulated using
complex Goursat functions, resulting in a set of singular integro-differential
equations on the boundary. The resulting boundary integral equations
are solved numerically using a Chebyshev collocation method which
is augmented by a fractional power term, derived by asymptotic analysis
of the corner region, to resolve stress singularities at corners of
the hole. We apply our numerical method to the test case of the hole
formed by two partially-overlapping circles, which can include either
a corner pointing into the solid or a corner pointing out of the solid.
Our numerical results recover the exact stress on the boundary to
within relative error $10^{-3}$ for modest computational effort.
\end{abstract}

\section{Introduction\label{sec:Introduction}}

Free boundary elasticity problems are fundamental to describing crystal
growth in strained solids. Due to the slow time scale of crystal growth
relative to the time scale of elastic relaxation, the elastic response
can be described by the quasi-static elasticity problem: a time-dependent
free or moving boundary problem for the morphology of the solid surface
coupled to the static elasticity equations in the solid. Thus, efficient
computational methods for solving the static elasticity problem for
general boundary shapes are necessary. 

The formulation of the elasticity problem for a given system with
a free or moving boundary can correspond to an interior domain, an
exterior domain, or a semi-infinite domain with boundaries that are
smooth or allowed to have corners. More complicated systems can have
multiple domains of elastically-interacting solid phases. Our focus
here is on the fundamental geometry of a hole (or void) inside a two-dimensional
elastically-stressed solid, for the case when the domain has corners
and the stress field has singularities at the corner.

There are of course exact solutions for the stressed infinite plate with
holes of different geometries. But many of these results apply only
to pre-defined hole geometries without corners (e.g. a circle \cite{Muskhelishvili}
or ellipse \cite{kolosoff1914some}), or for pre-defined geometries
in which there is some imposed geometric rounding of the corner (eg.
a rectangle with rounded corners \cite{savin1970stress,pan2013stress,motok1997stress}),
or in special hole-with-corner configurations (e.g. a hole formed
by partially overlapping circles \cite{ling1948stresses}). But all
of these results
are for holes of fixed geometry and thus
require the hole shape to be known in advance and
are not useful for solving a free boundary problem. In principle,
the elasticity solution for an arbitrarily shaped hole can be determined
by using a conformal map of the hole boundary to a circle (e.g. Schwarz\textendash Christoffel
mapping \cite{driscoll2002schwarz}), but in most cases the conformal
map approach will not work if the original domain has a corner. 

One approach to the problem of finding the solution for a domain with
corners is to impose some mechanism for rounding the corner which
gives a smooth boundary for the hole from which the elasticity solution
can be found. For example, the rectangular-hole solution of \cite{savin1970stress,pan2013stress,soutas2012elasticity01}
includes a geometrically-imposed corner-rounding radius. Corner-rounding
for a free-boundary problem can be naturally achieved by specifying
a curvature dependent surface energy that penalizes the formation
of corners and results in corner-rounding \cite{DiCarlo,golovin,gurtin1993thermomechanics}.
In fact, \cite{siegel2004evolution} uses this method to solve the
elastic stresses for a void of arbitrary shape and then combines the
elasticity solution to find the overall energy-minimizing void shape,
as well as dynamics for void shape changes due to mass transport.
While such corner-rounding methods permit solution of the elasticity
problem because the boundary of the solid is rendered smooth, the
corner-rounding also removes the weak (integrable) singularity of
the stress at the corner \cite{williams1952stress}. Of particular
interest to us is how the stress singularity of the corner does or
does not modify the behavior of the free boundary problem.\cite{chiu2020model,wu1982unconventional}

The role of a singularity in the elastic stress energy density near
a corner could in principle contribute to the energy balance determining
the free boundary equilibrium shape and potentially modify the equilibrium
corner angle. In the absence of elastic stress, the corner angle on
an equilibrium shape is given by a specific condition \cite{burton1951n,cabrera1964equilibrium}.
\cite{srolovitz2001stresses} consider the elastic energy locally
near the corner and use scaling arguments to argue that stresses do
not modify the corner angle from the no-stress results. In contrast,
\cite{siegel2004evolution} consider the energy-minimization problem
for the shape of a void in an elastic solid with anisotropic surface
energy by using a corner energy regularization term in which there
is an energy penalty for corners, and find that the apparent corner
angle does depend on elastic stress. So, to resolve the apparent discrepancy
in the influence of elastic stress on corners of energy-minimizing
free-boundary void shapes, we develop here a reliable numerical method
with high accuracy to determine the stress distribution of a void
with corners. These results are useful on their own, as a contribution
to understanding stress distribution due to voids with specific geometries,
and also as a necessary component in the more general problem of finding
the energy-minimizing void shape and understanding the fundamental
problem of the effect of elasticity on equilibrium corner angles.

The organization of the paper is as follows: in Section \ref{sec:Mathematical-formulation},
we derive a boundary integro-differential equation from the mathematical
formulation of the elasticity problem in an infinite plate under biaxial
stress; in Section 3 we describe the numerical method to discretize
the integro-differential equation in Section 2; in Section 4, we give
examples of the numerical results corresponding to hole shapes with and without
corners and analyze the error. Section 5 includes
a discussion of our numerical method and conclusion.

\section{Mathematical formulation \label{sec:Mathematical-formulation}}

\subsection{Boundary integro-differential equation}

We follow Muskhelishvili's complex variable formulation for two-dimensional
elasticity \cite{mikhlin1957integral,Muskhelishvili}. We consider
the exterior elasticity problem with a simply connected void and biaxial
stress applied at infinity. We assume plane-strain elasticity in the
$xy$-plane. The displacement field in the elastic solid is $\boldsymbol{u}(x,y)=u_{1}(\text{\ensuremath{x}},\text{\ensuremath{y}})\boldsymbol{e_{1}}+\text{\ensuremath{u}}_{2}(\text{\ensuremath{x}},\text{\ensuremath{y}})\boldsymbol{e_{2}}$,
where $\boldsymbol{e_{1}}$ and $\boldsymbol{e_{2}}$ are unit vectors
in  the $x,y$ directions. Let $D$ denote the solid region, $D'$
denote the void region, and let $\partial D$ represent the interface
between void and solid. Displacements and the stress tensor are defined
on $D$ and $\partial D$.

The infinitesimal strain tensor is defined by $E=\frac{1}{2}(\nabla\boldsymbol{u}+\nabla\boldsymbol{u}^{T})$,
and the first Piola-Kirchhoff stress tensor for the linearly elastic solid is given by 
$\text{\ensuremath{P}}=\lambda(\mbox{tr}(E))I+\mu(E+E^{T})$,
where $\mbox{tr}(E)$ is trace of matrix $E$,
$I$ is the identity tensor, and where $\lambda$
and $\mu$ are the Lame coefficients. Thus,
\begin{equation}
P=\left(\begin{array}{cc}
\sigma_{x} & \tau_{xy}\\
\tau_{xy} & \sigma_{y}
\end{array}\right)
\end{equation}
in the $xy$-plane with far-field condition:
\begin{equation}
\left(\begin{array}{cc}
\sigma_{x} & \tau_{xy}\\
\tau_{xy} & \sigma_{y}
\end{array}\right)\rightarrow\left(\begin{array}{cc}
1 & 0\\
0 & \chi
\end{array}\right)\quad\mbox{as}\;\sqrt{x^{2}+y^{2}}\rightarrow\infty,\label{eq:2}
\end{equation}
where $\chi=\sigma_{1}/\sigma_{2}$ is a parameter after nondimensionalizing
the stress components by the $x$-component of the applied biaxial stress.
Mechanical equilibrium in the solid gives (here $x_{1},x_{2}$ are
$x,y$):
\begin{equation}
\nabla\cdot P=\mathop{\underset{i,j=1,2}{\sum}}\frac{\partial P_{ij}}{\partial x_{j}}\boldsymbol{e_{i}}=0,\label{eq:3}
\end{equation}
\begin{equation}
P\cdot\mathbf{n}=0\quad\mbox{on}\quad\partial D,
\end{equation}
where $\mathbf{n}$ is the unit normal vector exterior to solid. Introduce
the stress function $W(x,y)$ as a smooth function defined on $D$
and $\partial D$, such that $\sigma_{x}=\partial^{2}W/\partial y^{2}$,
$\tau_{xy}=-\partial^{2}W/\partial x\partial y$, $\sigma_{y}=\partial^{2}W/\partial x^{2}$
(see \cite{gonzalez2008first}). Then Eq.\,\eqref{eq:3} is satisfied.
The compatibility condition for strain dictates that $W(x,y)$ satisfies
the biharmornic equation:
\begin{equation}
\frac{\partial^{4}W}{\partial x^{4}}+2\frac{\partial^{4}W}{\partial x^{2}\partial y^{2}}+\frac{\partial^{4}W}{\partial y^{4}}=0.
\end{equation}
Then using two functions $\phi(z)$ and $\psi(z)$ (called Goursat
functions) which are holomorphic on $D$ and $\partial D$ to represent
$W(x,y)$ with complex variable $z=x+iy$, we let $\text{\ensuremath{W}}(\text{\ensuremath{x}},\text{\ensuremath{y}})=\text{\mbox{Re}}\left\{ \bar{z}\phi(z)+\varsigma(z)\right\} $,
and $\psi(z)=\varsigma'(z)$. The relations between the stress components and
Goursat functions are then
\begin{equation}
\sigma_{x}+\sigma_{y}=4\mbox{Re}\left\{ \phi'(z)\right\} ,\label{eq:5}
\end{equation}
\begin{equation}
\sigma_{y}-\sigma_{x}+2i\tau_{xy}=2\left[\bar{z}\phi''(z)+\psi'(z)\right].\label{eq:6}
\end{equation}
Substituting Eq.\,\eqref{eq:5},\,\eqref{eq:6} into Eq.\,\eqref{eq:2},
the boundary conditions at infinity in terms of $\phi$ and $\psi$
become
\begin{equation}
\phi(z)=\frac{1+\chi}{4}z+C_{1}zi+C_{2}+O(\frac{1}{z})\quad\mbox{as}\;\mid z\mid\rightarrow\infty,\label{eq:7}
\end{equation}
\begin{equation}
\psi(z)=\frac{\chi-1}{2}z+C_{3}+O(\frac{1}{z})\quad\mbox{as}\;\mid z\mid\rightarrow\infty,\label{eq:8}
\end{equation}
where $C_{1}$ is arbitrary real constant and $C_{2},C_{3}$ are arbitrary
complex constants. To make the solution unique and for convenience,
we choose these arbitrary constants in $\phi$ and $\psi$ to be
zero, which does not affect the stresses (see \cite{Muskhelishvili,siegel2004evolution}).
Since no external force is applied on $\partial D$, the boundary condition
on $\partial D$ is given by
\begin{equation}
\phi(z)+z\overline{\phi'(z)}+\overline{\psi(z)}=0\quad\mbox{\mbox{on}}\;z\in\partial D.\label{eq:9}
\end{equation}
The Goursat functions $\phi$, $\psi$ can be written as $\phi(z)=(1+\chi)z/4+\varphi(z)$
and $\psi(z)=(\chi-1)z/2+h(z)$. Substitute into Eq.\,\eqref{eq:9}
and take the conjugate on both sides, then far-field conditions Eq.\,\eqref{eq:7},\,\eqref{eq:8}
and boundary condition Eq.\,\eqref{eq:9} are equivalent to
\begin{equation}
\varphi(z)\rightarrow0\quad\mbox{as}\;\mid z\mid\rightarrow\infty,\label{eq:10}
\end{equation}
\begin{equation}
h(z)\rightarrow0\quad\mbox{as}\;\mid z\mid\rightarrow\infty,\label{eq:11}
\end{equation}
\begin{equation}
\overline{\varphi(z)}+\frac{1+\chi}{4}\overline{z}+\overline{z}(\varphi'(z)+\frac{1+\chi}{4})+h(z)
+\frac{\chi-1}{2}z=0\quad\mbox{on}\;z\in\partial D.\label{eq:12}
\end{equation}
The purpose of making the substitution is to remove the singularity
of $\phi$ and $\psi$ at $\infty$. Then $\varphi$ and $h$ are analytic
on the region $D\cup\infty$. Multiply both sides of Eq.\,\eqref{eq:12}
by the factor $1/2\pi i\cdot dz/(z-t)$, where $t$ is an arbitrary
point in $D'$, and integrate along boundary $\partial D$, denoting
the integration contour $L$ as $\partial D$ traversed in the counterclockwise
direction. Since $\varphi$ and $h$ are analytic on the region $D\cup\infty$
with the conditions at $\infty$ (Eq.\,\eqref{eq:10},\,\eqref{eq:11}),
and $z$ is analytic in $D'$, by the Cauchy integral formula, the value of the resulting
Cauchy integrals are given by
\[
\frac{1}{2\pi i}\underset{L}{\int}\frac{\varphi(z)}{z-t}\,dz=0, \qquad\frac{1}{2\pi i}\underset{L}{\int}\frac{h(z)}{z-t}\,dz=0,
 \qquad\frac{1}{2\pi i}\underset{L}{\int}\frac{z}{z-t}\,dz=t.
\]
Eq.\,\eqref{eq:12} then becomes an integral equation which does not involve
$h(z)$:
\begin{equation}
\frac{1}{2\pi i}\underset{L}{\int}\frac{\overline{\varphi(z)}}{z-t}\,dz+\frac{1+\chi}{2\pi i\cdot2}\underset{L}{\int}\frac{\overline{z}}{z-t}\,dz
+\frac{1}{2\pi i}\underset{L}{\int}\frac{\overline{z}\varphi'(z)}{z-t}\,dz+\frac{\chi-1}{2}t=0.
\end{equation}
Now letting $t\rightarrow z_{0}$, where $z_{0}$ is a point on boundary
$\partial D$, the limits of the Cauchy integrals have the following properties:
\begin{equation}
\underset{t\rightarrow z_{0}}{\lim}\frac{1}{2\pi i}\underset{L}{\int}\frac{\overline{\varphi(z)}}{z-t}\,dz=\frac{1}{2}\overline{\varphi(z_{0})}+\frac{1}{2\pi i}\underset{L}{\int}\frac{\overline{\varphi(z)}}{z-z_{0}}\,dz,\label{eq:15}
\end{equation}

\begin{equation}
\underset{t\rightarrow z_{0}}{\lim}\frac{1}{2\pi i}\underset{L}{\int}\frac{\overline{z}}{z-t}\,dz=\frac{1}{2}\overline{z_{0}}+\frac{1}{2\pi i}\underset{L}{\int}\frac{\overline{z}}{z-z_{0}}\,dz,\label{eq:16}
\end{equation}

\noindent 
\begin{equation}
\underset{t\rightarrow z_{0}}{\lim}\frac{1}{2\pi i}\underset{L}{\int}\frac{\overline{z}\varphi'(z)}{z-t}\,dz=\frac{1}{2}\overline{z_{0}}\varphi'(z_{0})+\frac{1}{2\pi i}\underset{L}{\int}\frac{\overline{z}\varphi'(z)}{z-z_{0}}\,dz.\label{eq:17}
\end{equation}
All integrals on the right side of Eq.\,\eqref{eq:15}-\eqref{eq:17}
are in the sense of Cauchy principle value. Thus we obtain a singular
boundary integro-differential equation on $\partial D$:
\begin{multline}
\frac{1}{2}\overline{\varphi(z_{0})}+\frac{1}{2\pi i}\underset{L}{\int}\frac{\overline{\varphi(z)}}{z-z_{0}}\,dz+\frac{1+\chi}{4}\overline{z_{0}}+\frac{1+\chi}{4\pi i}\underset{L}{\int}\frac{\overline{z}}{z-z_{0}}\,dz\\
+\frac{1}{2}\overline{z_{0}}\varphi'(z_{0})+\frac{1}{2\pi i}\underset{L}{\int}\frac{\overline{z}\varphi'(z)}{z-z_{0}}\,dz+\frac{\chi-1}{2}z_{0}=0.\label{eq:18}
\end{multline}
Once the integro-differential equation is solved for $\varphi(z)$ on
the boundary $\partial D$, Eq.\,\eqref{eq:12} determines $h$ on the boundary.
The stress at any point $\zeta$ inside solid can then be determined
by analytic continuation of 
boundary values $\varphi(z)$ and $h(z)$ into
the domain $D$ using
\begin{equation}
\varphi(\zeta)=\frac{1}{2\pi i}\underset{-L}{\int}\frac{\varphi(z)}{z-\zeta}dz, \quad h(\zeta)=\frac{1}{2\pi i}\underset{-L}{\int}\frac{h(z)}{z-\zeta}dz.
\end{equation}
 The analytic continuation of $\varphi$ and $h$ into the domain can then
be used to construct the stress tensor $P$ in the solid from Eq.\,\eqref{eq:5},\,\eqref{eq:6}.

\subsection{Local asymptotic analysis near corners\label{subsec:2.2}}

\begin{figure}
\centering
\includegraphics[width=8.4cm]{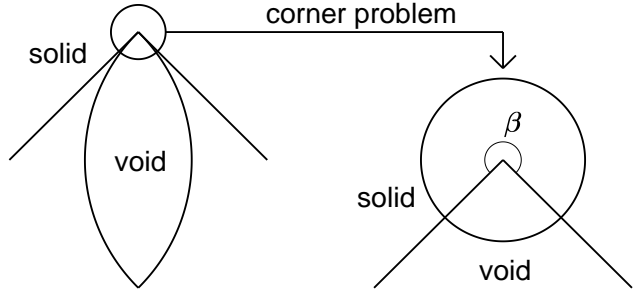}

\caption{\label{fig:1}Boundary shape of the corner problem.}
\end{figure}

Goursat function $\varphi(z)$ is smooth along the boundary if the
boundary shape of the hole is smooth (with no corners). For a convex void
shape with a corner,  $\varphi'(z)$ has a singularity at the corner and the
stress goes to infinity when $z$ approaches the corner \cite{savin1970stress,williams1952stress}.
We analyze the stress asymptotically to determine the order of the
singularity as a function of corner angle. By setting the vertex of
the corner as the origin, the shape of the boundary of the hole near the corner is a wedge with same angle as the corner angle (see Fig.\,\ref{fig:1}).
The biharmornic equation in polar coordinates is given in \cite{soutas2012elasticity01}
as
\begin{equation}
(\frac{\partial^{2}}{\partial r^{2}}+\frac{1}{r}\frac{\partial}{\partial r}+\frac{1}{r^{2}}\frac{\partial^{2}}{\partial\theta^{2}})^{2}W=0,\label{eq:20}
\end{equation}
and the corresponding stress components in polar coordinates are given
by
\begin{equation}
\sigma_{rr}=\frac{1}{r}\frac{\partial W}{\partial r}+\frac{1}{r^{2}}\frac{\partial^{2}W}{\partial\theta^{2}},
\end{equation}
\begin{equation}
\sigma_{\theta\theta}=\frac{\partial^{2}W}{\partial r^{2}},
\end{equation}
\begin{equation}
\sigma_{r\theta}=-\frac{\partial}{\partial r}(\frac{1}{r}\frac{\partial W}{\partial\theta}),
\end{equation}
where $\sigma_{rr}$ is the stress in the radial direction, $\sigma_{\theta\theta}$
is the stress in the $\theta$ direction, and $\sigma_{r\theta}$ is the
shear stress.

Letting $r=\varepsilon\tilde{r}$ and $W(r,\theta)=\widetilde{W}(\tilde{r},\theta)$
with $\varepsilon\ll1$ to find the corner solution, the biharmonic
Eq.\,\eqref{eq:20} becomes:
\begin{equation}
\frac{1}{\varepsilon^{2}}(\frac{\partial^{2}}{\partial\tilde{r}^{2}}+\frac{1}{\tilde{r}}\frac{\partial}{\partial\tilde{r}}+\frac{1}{\tilde{r}^{2}}\frac{\partial^{2}}{\partial\theta^{2}})^{2}\widetilde{W}=0.\label{eq:24}
\end{equation}
Stresses of the corner problem are:
\begin{equation}
\widetilde{\sigma_{\tilde{r}\tilde{r}}}=\frac{1}{\tilde{r}}\frac{\partial\widetilde{W}}{\partial\tilde{r}}+\frac{1}{\tilde{r}^{2}}\frac{\partial^{2}\widetilde{W}}{\partial\theta^{2}}=\varepsilon^{2}\sigma_{rr},
\end{equation}
\begin{equation}
\widetilde{\sigma_{\theta\theta}}=\frac{\partial^{2}\widetilde{W}}{\partial\tilde{r}^{2}}=\varepsilon^{2}\sigma_{\theta\theta},
\end{equation}
\begin{equation}
\widetilde{\sigma_{\tilde{r}\theta}}=-\frac{\partial}{\partial\tilde{r}}(\frac{1}{\tilde{r}}\frac{\partial\widetilde{W}}{\partial\theta})=\varepsilon^{2}\sigma_{r\theta},
\end{equation}
The differential equation of the corner problem is still a biharmonic
equation from Eq.\,\eqref{eq:24}.  Since the far-field conditions $\sigma_{rr},\sigma_{r\theta},\sigma_{\theta\theta}$
are finite away from the corner, $\widetilde{\sigma_{rr}},\widetilde{\sigma_{r\theta}},\widetilde{\sigma_{\theta\theta}}\rightarrow0$
as $\tilde{r}\rightarrow\infty$. 
Hence the local corner problem is identical to the wedge problem.   The solution
to the wedge problem with ``free-free'' boundary conditions is given by separation of variables in \cite{williams1952stress}. The
order of $\varphi$, $\varphi'$ and stresses near the corner are
\begin{equation}
\sigma_{rr},\sigma_{r\theta},\sigma_{\theta\theta}\sim r^{\lambda-2},\varphi'(r,\theta)\sim r^{\lambda-2},\varphi(r,\theta)\sim r^{\lambda-1},\label{eq:28}
\end{equation}
where $\lambda$ is the solution of 
\begin{equation}
\sin\left[(\lambda-1)\beta\right]=-(\lambda-1)\sin\beta \label{eq:29}
\end{equation}
and where $\beta$ is the corner angle (see Fig.~\ref{fig:1}).
From Eq.\,\eqref{eq:29}, we can determine the behavior of the stresses
near the corner. In case of $\beta<\pi$, $\lambda$ is greater than
$2$ and there is no singularity near the corner. In case of $\beta=\pi$,
the boundary is a straight line and the stresses are constant as $\lambda=2$.
When $\pi<\beta<2\pi$ (see Fig.\,\ref{fig:1}), $\lambda$ is between
$1$ and $2$ which gives stress singularity near the corner. If $\beta$
goes to $2\pi$, $\lambda$ approaches 1. For all cases in which $\pi<\beta<2\pi$,
the singularity in the stresses near the corner is  an integrable singularity \cite{williams1952stress}.

\section{Numerical method}

\noindent  The surface of the elastic solid is described by a closed
continuous curve on the $xy$-plane. To illustrate and test our numerical
method we let the center of the void be the origin and consider two-fold-symmetric void shapes with reflection
symmetry across the $x$ and $y$ axes. 
We consider a class of shapes that are piecewise-smooth except for
possible corners on the $x$ and $y$ axes (see Fig.\,\ref{fig:8} for
an illustration). Such corners could be present in the full energy
minimization problem for the free boundary shape when the surface
energy model is strongly anisotropic with excluded orientations at
the corner corresponding to negative surface stiffness \cite{herring1951some}.
The assumed symmetry
of the void shape results in symmetry of the Goursat functions in
the complex plane. Thus we consider the shape and the Goursat functions
on the first quadrant and use symmetry to extend to the entire $xy$-plane.

\subsection{Modeling the surface\label{subsec:3.1}}

The surface is represented in polar coordinates $r(\theta)$, where
$r$ is the radial coordinate and $\theta\in[0,\pi/2]$ denotes the
polar angle (note that the polar coordinate is different from Section \ref{subsec:2.2}). 
For example, $r(\theta)=1$ when surface is a circle
or $r(\theta)=b/\sqrt{1-(e\cos\theta)^{2}}$ for an ellipse with eccentricity
$e$ and semi-minor axis $b$. If the surface is given analytically,
$r'(\theta)$ and $r''(\theta)$ can be obtained directly without
error. If the surface is not given explicitly, we use a series in Chebyshev
polynomials to represent the shape as 
\begin{equation}
r(\theta)=\stackrel[k=0]{N-1}{\sum}c_{k}T_{k}(\theta),
\end{equation}
where $T_{k}(\theta)$ is the $k$-th Chebyshev polynomial on $[0,\pi/2]$.
Coefficients $c_{k}$ are derived from the value of $r$ at Chebyshev
nodes on $\theta\in[0,\pi/2]$ using Chebyshev interpolation \cite{press1992numerical}.
For cases in which the surface
shape is smooth the error can be reduced to machine roundoff error
for sufficiently large N. 
More details about error analysis of Chebyshev interpolation are given in \cite{boyd2001chebyshev}. 

The problem for the elastic stress distribution has been reduced to
an integro-differential equation given by Eq.\,\eqref{eq:18} for
the boundary values of the complex function $\varphi(z)$ where $z(\theta)=r(\theta)\cos\theta+i\cdot r(\theta)\sin\theta$.
Goursat function $\varphi$ is holomorphic on $D$, which leads to
$\varphi$ being a continuous function of $\theta$ on the boundary.
We use a Chebyshev basis to represent the real and imaginary parts
of the boundary values of the Goursat function on the first quadrant:
\begin{equation}
\varphi(\theta)=\stackrel[k=0]{N-1}{\sum}a_{k}T_{k}(\theta)+i\cdot\stackrel[k=0]{N-1}{\sum}b_{k}T_{k}(\theta),\label{eq:31}
\end{equation}
where $a_{k}$, $b_{k}$ are unknowns. If the shape of the hole has
corners, from the local asymptotic analysis near corners in Section
2.2, the Goursat function $\varphi$ near the corner can not be well
approximated only by polynomials. Thus, considering the example where
the corner is located at $\theta_{c}=\pi/2$, we add a corner term
in $s^{\lambda-1}$, where $s$ is the arclength from the corner
and $\lambda$ is the solution of Eq.\,\eqref{eq:29}. We thus have
$\varphi(\theta)\sim c_{1}s^{\lambda-1}\sim c_{2}|\theta_{c}-\theta|^{\lambda-1}$
near the corner where $c_{1},c_{2}$ are constants.  To accommodate the 
local behavior near the corner we modify the expansion in Eq.\,\eqref{eq:31} as
\begin{equation}
\varphi(\theta)=a_{N-1}(\frac{\pi}{2}-\theta)^{\lambda-1}+\stackrel[k=0]{N-2}{\sum}a_{k}T_{k}(\theta)
+i\cdot\left(b_{N-1}(\frac{\pi}{2}-\theta)^{\lambda-1}+\stackrel[k=0]{N-2}{\sum}b_{k}T_{k}(\theta)\right).\label{eq:32}
\end{equation}
\noindent Using the symmetry of the void shape the associated symmetries
of the elastic stress are
\[
\sigma_{x}(x,y)=\sigma_{x}(-x,y)=\sigma_{x}(x,-y)=\sigma_{x}(-x,-y),
\]
\[
\sigma_{y}(x,y)=\sigma_{y}(-x,y)=\sigma_{y}(x,-y)=\sigma_{y}(-x,-y),
\]
\[
\tau_{xy}(x,y)=-\tau_{xy}(-x,y)=\tau_{xy}(-x,-y)=-\tau_{xy}(x,-y).
\]
Thus, from the dependence of the stress on the Goursat functions we
have that $\mbox{Re}\{\varphi(\theta)\}$ is odd in $x$ and even
in $y$ while $\mbox{Im}\{\varphi(\theta)\}$ is even in $x$ and
odd in $y$. We can use the symmetry to extend $\varphi(\theta)$
on $0<\theta<\pi/2$ to $\pi/2<\theta<2\pi$ as
\[
\varphi(\theta)=-\mbox{Re}\left\{ \varphi(\pi-\theta)\right\} +i\cdot\mbox{Im}\left\{ \varphi(\pi-\theta)\right\} \quad\mbox{on}\;\theta\in[\pi/2,\pi].
\]
\[
\varphi(\theta)=-\mbox{Re}\left\{ \varphi(\theta-\pi)\right\} -i\cdot\mbox{Im}\left\{ \varphi(\theta-\pi)\right\} \quad\mbox{on}\;\theta\in[\pi,3\pi/2].
\]
\[
\varphi(\theta)=\mbox{Re}\left\{ \varphi(2\pi-\theta)\right\} -i\cdot\mbox{Im}\left\{ \varphi(2\pi-\theta)\right\} \quad\mbox{on}\;\theta\in[3\pi/2,2\pi].
\]
Since $\varphi$ is a continuous function of $\theta$, we have the
continuity conditions on both ends of the interval
\begin{equation}
\mbox{Re}\left\{ \varphi(\theta)\right\} =0\quad\mbox{at}\;\theta=\pi/2,\label{eq:33}
\end{equation}
\begin{equation}
\mbox{Im}\left\{ \varphi(\theta)\right\} =0\quad\mbox{at}\;\theta=0\text{.}\label{eq:34}
\end{equation}
Note that the integro-differential equation \eqref{eq:18} admits a homogeneous
solution $\varphi_{H}(z)=iaz+b$, where $a$ is a real constant and
$b$ is a complex constant, which corresponds to arbitrary degrees
of freedom in the representation of the stress-free state by Goursat
functions \cite{mikhlin1957integral}. However, because of the assumed
symmetry of the shape and our resulting symmetry relations for the
real and imaginary parts of $\varphi(z)$, the symmetry excludes the
homogeneous solution $\varphi_{H}(z)$ because $b$ has only even
symmetry in both $x$ and $y$, while $iaz$ has only odd symmetry
for $x$ and $y$. 

Some more constraints will be discussed in Section \ref{subsec:3.4}
to ensure that $\varphi(\theta)$ is boundary value of an analytic
function.

\subsection{Nested Gauss-Legendre quadrature\label{subsec:3.2}}

Traditional Gauss\textendash Legendre quadrature \cite{golub1969calculation}
gives an approximation to the integral of function $f(x)$ on the interval
$[-1,1]$ as
\begin{equation}
\int_{-1}^{1}f(x)\,dx=\stackrel[i=1]{N}{\sum}\omega_{i}f(x_{i}),
\end{equation}
where $x_{i}$ are Gauss-Legendre quadrature points and $\omega_{i}$
are corresponding Gauss-Legendre quadrature weights. Gauss-Legendre
quadrature is exact for polynomials under degree $2N$  and it converges
as $N\rightarrow\infty$ for smooth functions $f(x)$ which can be
approximated by polynomials. However, in the corner case, Eq.\,\eqref{eq:32}
has the corner term with a non-integer exponent $\lambda-1$ when
the corner angle of the hole $\beta\neq\pi$. We introduce nested
Gauss-Legendre quadrature \cite{bremer2010nonlinear,hoskins2019numerical}
to improve the convergence of Gauss-Legendre quadrature for the corner
term. 

As an illustration of the application of nested Gauss-Legendre quadrature
we evaluate $\int_{0}^{1}f(x)dx$ with possible integrable singularity
near $x=1$. The nested Gauss-Legendre quadrature algorithm is \cite{bremer2010nonlinear,hoskins2019numerical}:\\
\textit{Step 1}: Find Gauss-Legendre quadrature points $x_{i}$ and
Gauss-Legendre quadrature weights $\omega_{i}$ on $[0,1]$ with $N$
quadrature points. Set $a_{1}=0$ and number of iterations $n=1$.\\
\textit{Step }2: Evaluate $\int_{a_{n}}^{1}f(x)dx$ using Gauss-Legendre
quadrature on $[a_{n},1]$.\\
\textit{Step }3: Let $a_{n+1}=(1+a_{n})/2$. Divide $[a_{n},1]$ into
two sub-intervals $[a_{n},a_{n+1}]$ and $[a_{n+1},1]$.\\
\textit{Step }4: Evaluate $\int_{a_{n}}^{a_{n+1}}f(x)\,dx$ and $\int_{a_{n+1}}^{1}f(x)\,dx$
using Gauss-Legendre quadrature.\\
\textit{Step }5: Repeat \textsl{Step 2-4} until 
\[
\left|\int_{a_{n}}^{1}f(x)\,dx-\left[\int_{a_{n}}^{a_{n+1}}f(x)\,dx+\int_{a_{n+1}}^{1}f(x)\,dx\right]\right|<\varepsilon,
\]
where $\varepsilon$ is a prescribed tolerance. Then 
\[
\int_{0}^{1}f(x)\,dx=\stackrel[i=1]{n}{\sum}\int_{a_{i}}^{a_{i+1}}f(x)\,dx+\int_{a_{n+1}}^{1}f(x)\,dx.
\]

\begin{figure}
\centering
\includegraphics[width=8.4cm]{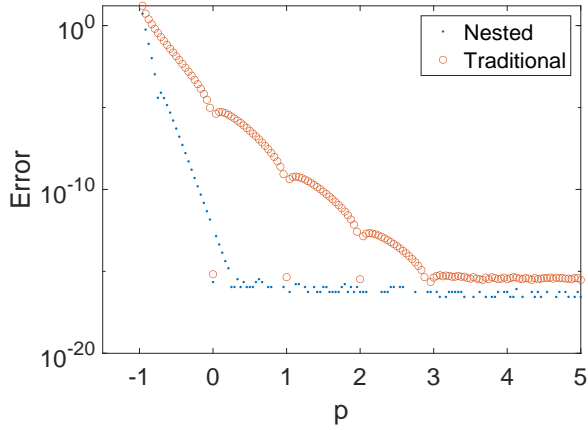}

\caption{\label{fig:2}Comparison of traditional Guass-Legendre quadrature
($N=16$) and nested Gauss-Legendre quadrature ($N=16$ and $\varepsilon=10^{-15}$)
for the integral of $x^{p}$ on $(0,1)$. }
\end{figure}

As a specific test case, we evaluate $\int_{0}^{1}(1-x)^{p}dx=1/(p+1)$
with $-1<p\leq5$ using both traditional and nested Gauss-Legendre
quadrature. Fig.\,\ref{fig:2} shows the semi-log plot of the error
versus $p$. Traditional Gauss-Legendre quadrature uses a polynomial
basis on the given interval to evaluate the integral. Then the quadrature
is exact for functions on the polynomial space. When $0<p<3$ and
$p$ is not integer, traditional Gauss-Legendre quadrature has significant
error. Nested Gauss-Legendre quadrature uses a piecewise polynomial
to approximate the function near end of the interval and has better
performance in this case. In the singular case ($-1<p<0$), both methods
have significant error but the nested Gaussian quadrature gives substantially
smaller errors. Since Eq.\,\eqref{eq:32} has corner terms with non-integer
power, we use nested Gauss-Legendre quadrature to evaluate the integrals.

\subsection{Discretization of integral equations\label{subsec:3.3}}

Eq.\,\eqref{eq:18} holds for any $z_{0}$ on $L$. We pick $z_{0}=z_{i}=z(\theta_{i})$
$i=1,2,\ldots,N-1$ as collocation points, where $\theta_{i}=\pi(x_{i}+1)/4$
and $x_{i}$ is the i-th root of degree ($N-1$) Legendre polynomial
$P_{N-1}$ (Gauss-Legendre quadrature points). The choice of the number
of collocation points determines the number of unknowns. We implement
both real part and imaginary part of Eq.\,\eqref{eq:18} at each $z_{i}$
($2(N-1)$ equations) together with two boundary conditions \eqref{eq:33} and
\eqref{eq:34} to give a system of $2N$ equations for a total of $2N$
unknowns $a_{k}$, $b_{k}$. 
The integrals in Eq.\,\eqref{eq:18}
are evaluated numerically using nested Gauss-Legendre quadrature
described in Section \ref{subsec:3.2}. $\varphi(\theta_{i})$ can
be obtained from the Chebyshev coefficients by using Clenshaw's recurrence formula from \cite{clenshaw1955note}
to minimize the truncation error when evaluating the Chebyshev series
at the collocation points $\theta=\theta_{i}$. 

Now we consider the integrals in Eq.\,\eqref{eq:18} term by term.
The singularity of the Cauchy principal value integrals in Eq.\,\eqref{eq:18}
is extracted using:
\begin{equation}
\underset{L}{\int}\frac{\overline{\varphi(z)}}{z-z_{i}}\,dz=\underset{L}{\int}\frac{\overline{\varphi(z_{i})}}{z-z_{i}}\,dz+\underset{L}{\int}\frac{\overline{\varphi(z)}-\overline{\varphi(z_{i})}}{z-z_{i}}\,dz.\label{eq:36}
\end{equation}
The first term is
\begin{equation}
\underset{L}{\int}\frac{\overline{\varphi(z_{i})}}{z-z_{i}}\,dz=\overline{\varphi(z_{i})}\cdot\underset{L}{\int}\frac{1}{z-z_{i}}\,dz=\pi i\thinspace\overline{\varphi(z_{i})}
\end{equation}
from boundary version of Cauchy integral formula. The second term
takes the form:
\begin{equation}
\underset{L}{\int}\frac{\overline{\varphi(z)-\varphi(z_{i})}}{z-z_{i}}\,dz=\int_{0}^{2\pi}\frac{\overline{\varphi(\theta)-\varphi(\theta_{i})}}{z(\theta)-z(\theta_{i})}\cdot z'(\theta)\,d\theta,
\end{equation}
where 
\[
z'(\theta)=r'(\theta)\cos\theta-r(\theta)\sin\theta+i(r'(\theta)\sin\theta+r(\theta)\cos\theta).
\]
Notice the singularity at $\theta=\theta_{i}$ is a removable singularity
canceled by taking limit at $\theta_{i}$
\begin{equation}
\underset{\theta\rightarrow\theta_{i}}{\lim}\frac{\overline{\varphi(\theta)-\varphi(\theta_{i})}}{z(\theta)-z(\theta_{i})}=\frac{\overline{\varphi'(\theta_{i})}}{z'(\theta_{i})}.
\end{equation}
We find $\varphi'(\theta_{i})$ by using the algorithm in \cite{press1992numerical}
based on the relation between the Chebyshev coefficients of a function
and the Chebyshev coefficients of its derivatives as
\begin{equation}
c'_{i-1}=c'_{i+1}+2(i-1)c_{i}\quad i\geq1,\label{eq:40}
\end{equation}
where $c_{i}$ are the Chebyshev coefficients and $c'_{i}$ are Chebyshev
coefficient of the derivatives. The Chebyshev coefficients of $\varphi'$
are linear in unknowns $a_{k}$, $b_{k}$ from Eq.\,\eqref{eq:40}.
Then integral \eqref{eq:36} can be evaluated using nested Gauss-Legendre
quadrature and implementing symmetry of $\varphi(\theta)$ and $r(\theta)$.
Continuing with the next integral term in Eq.\,\eqref{eq:18} we evaluate the $\overline{z}$ integral in Eq.\,\eqref{eq:40} the
same way, giving
\begin{equation}
\underset{L}{\int}\frac{\overline{z}}{z-z_{i}}\,dz=\pi i\thinspace\overline{z_{i}}+\int_{0}^{2\pi}\frac{\overline{z}-\overline{z_{i}}}{z-z_{i}}z'(\theta)\,d\theta.
\end{equation}
Next we consider the $\overline{z}\varphi'(z)$ term in Eq.\,\eqref{eq:18}:
\begin{equation}
\underset{L}{\int}\frac{\overline{z}\varphi'(z)}{z-z_{i}}\,dz=\underset{L}{\int}\frac{\overline{z_{i}}\varphi'(z_{i})}{z-z_{i}}\,dz+\underset{L}{\int}\frac{\overline{z}\varphi'(z)-\overline{z_{i}}\varphi'(z_{i})}{z-z_{i}}\,dz.
\end{equation}
Following the same logic, the first term is 
\begin{equation}
\underset{L}{\int}\frac{\overline{z_{i}}\varphi'(z_{i})}{z-z_{i}}\,dz=\pi i\thinspace\overline{z_{i}}\varphi'(z_{i}).
\end{equation}
Rewriting the second integral in terms of $\theta$ gives
\begin{equation}
\underset{L}{\int}\frac{\overline{z}\varphi'(z)-\overline{z_{i}}\varphi'(z_{i})}{z-z_{i}}\,dz
=\int_{0}^{2\pi}\frac{\overline{z(\theta)}\varphi'(\theta)-\overline{z(\theta_{i})}\varphi'(\theta_{i})}{z(\theta)-z(\theta_{i})}\,d\theta.\label{eq:44}
\end{equation}
The singularity at $\theta=\theta_{i}$ is a removable singularity
with
\begin{equation}
\underset{\theta\rightarrow\theta_{i}}{\lim}\frac{\overline{z(\theta)}\varphi'(\theta)-\overline{z(\theta_{i})}\varphi'(\theta_{i})}{z(\theta)-z(\theta_{i})}=\frac{\overline{z'(\theta_{i})}\varphi'(\theta_{i})+\overline{z(\theta_{i})}\varphi''(\theta_{i})}{z'(\theta_{i})}.
\end{equation}
Nested Gauss-Legendre quadrature has significant error when function
has a singularity. Therefore, Eq.\,\eqref{eq:44} needs to be evaluated
carefully near the corner using integration by parts to avoid the
singularity of $\varphi'(\theta)$ at $\theta=\pi/2$:

\begin{equation}
\int_{0}^{\pi/2}\frac{\overline{z(\theta)}\varphi'(\theta)}{z(\theta)-z(\theta_{i})}\,d\theta=\int_{0}^{\pi/2-\theta_{\varepsilon}}\frac{\overline{z(\theta)}\varphi'(\theta)}{z(\theta)-z(\theta_{i})}\,d\theta
+\int_{\pi/2-\theta_{\varepsilon}}^{\pi/2}\frac{\overline{z(\theta)}\varphi'(\theta)}{z(\theta)-z(\theta_{i})}\,d\theta,
\end{equation}
\begin{equation}
\int_{\pi/2-\theta_{\varepsilon}}^{\pi/2}\frac{\overline{z(\theta)}\varphi'(\theta)}{z(\theta)-z(\theta_{i})}\,d\theta=\frac{\overline{z(\theta)}\varphi(\theta)}{z(\theta)-z(\theta_{i})}\biggr\rvert_{\pi/2-\theta_{\varepsilon}}^{\pi/2}
-\int_{\pi/2-\theta_{\varepsilon}}^{\pi/2}\varphi(\theta)\,d\left(\frac{\overline{z(\theta)}}{z(\theta)-z(\theta_{i})}\right),
\end{equation}
where $\theta_{\varepsilon}=\theta_{1}/2$, and $\theta_{1}$ is the
polar angle at first collocation point. 

We substitute all results in this section to discretize Eq.\,\eqref{eq:18}.
$\varphi(\theta)$ is a linear combination of unknowns $a_{k}$ and
$b_{k}$ for any $\theta$ from Eq.\,\eqref{eq:32}. The coefficients
of derivative of a Chebyshev approximated function are linear in the
coefficients of original Chebyshev approximated function (Eq.\,\eqref{eq:40}),
thus $\varphi'(\theta_{i})$ and $\varphi''(\theta_{i})$ are linear
in $a_{k}$ and $b_{k}$. Gauss-Legendre quadrature is a weighted sum
of the function values at quadrature points, which is a linear operator
in $\varphi$ or $\varphi'$. Finally, we get a linear system of $a_{k}$
and $b_{k}$ from discretizing Eq.\,\eqref{eq:18}.

\subsection{Analyticity equations\label{subsec:3.4}}

In Section \ref{subsec:3.1}, we assume $\varphi$ is a smooth function
of $\theta$ on first quadrant of boundary $\partial D$ from \eqref{eq:31}.
However, $\varphi$ in Eq.\,\eqref{eq:18} is an analytic function
on $D\cup\partial D$. Not every function with smooth real and imaginary
part on the boundary is an analytic function. So more constraints are
needed to make $\varphi$ an analytic function in $D$. Since $\varphi(\zeta)$
is analytic for $\zeta$ in $D$, by the Cauchy integral formula and
our requirement that $\varphi(\infty)=0$ from Eq.\,\eqref{eq:10},
\begin{equation}
\varphi(\zeta)=\frac{1}{2\pi i}\underset{-L}{\int}\frac{\varphi(z)}{z-\zeta}\thinspace dz.
\end{equation}
Thus, 
\begin{equation}
\varphi(\zeta)=\frac{1}{2\pi i}\underset{-L}{\int}\frac{\varphi(z)}{z-\zeta}\,dz\rightarrow0\quad\mbox{as}\,\mid\zeta\mid\rightarrow\infty.\label{eq:49}
\end{equation}
Eq.\,\eqref{eq:49} holds because $\varphi(z)$ is bounded on $\partial D$.
This definition of $\varphi(\zeta)$ guarantees analyticity in $D$.
To make $\varphi(\zeta)$ analytic on $D\cup\partial D$, $\varphi$
should be continuous to any point on boundary $\partial D$. Since
collocation points $z_{i}$ are on boundary $\partial D$, $\varphi(\zeta)$
is continuous as $\zeta\rightarrow z_{i}$ for all $i$:
\begin{equation}
\underset{\zeta\rightarrow z_{i}}{\lim}\frac{1}{2\pi i}\underset{-L}{\int}\frac{\varphi(z)}{z-\zeta}\,dz=\varphi(z_{i}).
\end{equation}
Take the limit and evaluate the integral using its Cauchy principle
value, the limit becomes: 
\begin{equation}
\frac{1}{2\pi i}\underset{L}{\int}\frac{\varphi(z)}{z-z_{i}}\thinspace dz+\frac{1}{2}\varphi(z_{i})=0\label{eq:51}
\end{equation}
for any collocation point $z_{i}$. The integral in Eq.\,\eqref{eq:51}
can be evaluated in the same way as we find the Cauchy integral of
$\overline{\varphi(z)}$ in section 3.2. Then Eq.\,\eqref{eq:51}
at collocation points are a set of linear equations of $a_{k}$ and
$b_{k}$. Since Eq.\,\eqref{eq:51} on all points on $\partial D$ enforces
the analyticity of $\varphi(\zeta)$ on $D\cup\partial D$,  we therefore
use Eq.\,\eqref{eq:51} at each collocation point in addition to
our boundary integral equation \eqref{eq:18} to construct our analytic $\varphi(z)$.

We know that the solution to the two-dimensional elasticity problem
is unique \cite{knops2012uniqueness}. But if we solve Eq.\,\eqref{eq:18}
without constraints Eq.\,\eqref{eq:51}, other non-analytic $\varphi$
cause the solution to be non-unique. Eq.\,\eqref{eq:51} at collocation
points restricts the solution to the space of boundary values of analytic
functions, which ensures the uniqueness of the solution to discretized
problem. Finally, by combining the $2(N-1)$ equations from Eq.\,\eqref{eq:51}
at collocation points, the $2(N-1)$ equations from discretized integral
equation Eq.\,\eqref{eq:18}, and two equations Eq.\,\eqref{eq:33},\,\eqref{eq:34}
from boundary conditions at each end of the domain, we obtain an overdetermined
linear system of $(4N-2)$ equations for the $2N$ unknowns $a_{k}$ and
$b_{k}$. The non-square linear system can be solved in sense of least
squares using the MATLAB matrix left division operator (QR method).

\section{Numerical results}

Now we have a numerical method to determine the Goursat function on
the boundary. In this section, we test our numerical method on cases
for which the boundary shape corresponds to a circle, an ellipse,
and overlapping circles.

\subsection{Measurements of the error\label{subsec:4.1}}

In the cases of the circle and the ellipse, the exact solution for
the Goursat function on the boundary are given in \cite{Muskhelishvili}.
We can compare $\varphi$ from our numerical method with exact solution
using the $L^{2}$ error norm:
\begin{equation}
\mbox{Error}_{L^{2}}=\left[\frac{2}{\pi}\int_{0}^{\pi/2}\lvert\varphi(\theta)-\varphi_{\mbox{exact}}(\theta)\lvert{}^{2}\,d\theta\right]^{1/2}.
\end{equation}
In the case of overlapping circles, there is an exact solution for
the trace of the stress $\sigma_{x}+\sigma_{y}$ \cite{ling1948stresses}
but not for $\varphi$. Comparing $\sigma_{x}+\sigma_{y}=1+\chi+4\mbox{Re}\left\{ \varphi'(z)\right\} $
with the exact solution is a good numerical test for our method. We
thus define the $L^{2}$ norm error of $\sigma_{x}+\sigma_{y}$ as
\begin{equation}
\mbox{Error}_{L^{2}}=\left[\frac{2}{\pi}\int_{0}^{\pi/2}((\sigma_{x}+\sigma_{y})-(\sigma_{x}+\sigma_{y})_{\mbox{exact}}){}^{2}\,d\theta\right]^{1/2}.
\end{equation}

\subsection{Test for circle\label{subsec:4.2}}

\begin{figure}
\centering
\includegraphics[width=8.4cm]{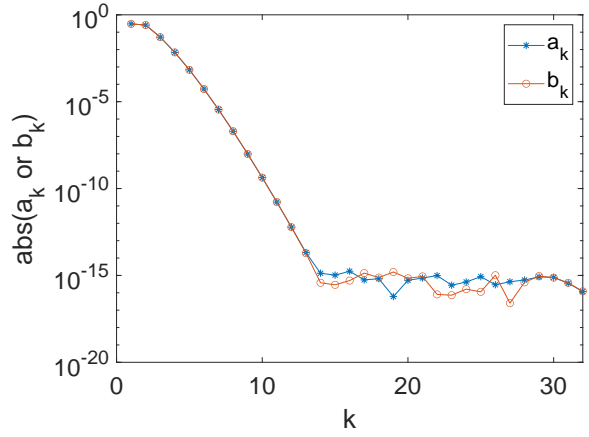}

\caption{\label{fig:3}Result for test case of a circle: absolute value of
coefficient of $a_{k}$ or $b_{k}$ versus index $k$ for the number
of collocation points $N=32$.{\footnotesize{} }}
\end{figure}

\begin{figure}
\centering
\includegraphics[width=8.4cm]{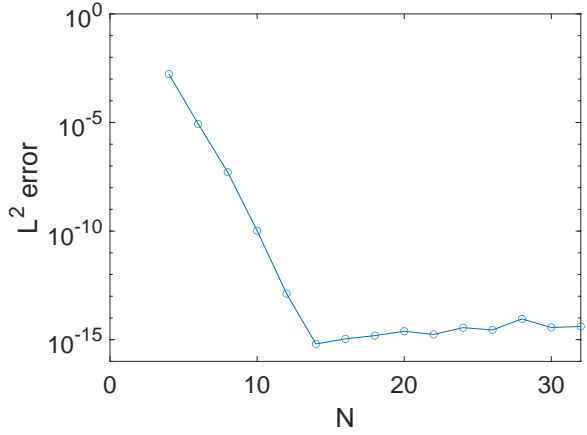}

\caption{\label{fig:4}Result for test case of a circle: $L^{2}$ error of
$\varphi(z)$ versus number of collocation points $N$. Here $L^{2}$
error is defined in Section \ref{subsec:4.1}.}
\end{figure}

When the interface is a unit circle, the exact solution of the interior
problem for the case of uniaxial tension, $\chi=0$, is $\varphi(z)=2z$.
The solution of the exterior problem is derived by the conformal map
$z(\zeta)=1/\zeta$ which maps the interior domain to the exterior
domain. Here $r(\theta)=1$ and $r'(\theta)=r''(\theta)=0$. The exact
solution of exterior problem for the case of uniaxal tension $\chi=0$
is given by $\varphi_{\mbox{exact}}(z)=1/(2z)$ \cite{Muskhelishvili}.
Absolute values of the computed Chebyshev coefficients for $\varphi(z)$,
$a_{k}$ and $b_{k}$, are presented in Fig.\,\ref{fig:3}. Coefficients
of $\varphi(z)$ decrease to magnitudes comparable to the truncation
error at $N=15$ which confirms the convergence of our collocation method.
The convergence rate is shown in the $L^{2}$ error versus $N$ plot
in Fig.\,\ref{fig:4}. These results verify that our numerical algorithm
has spectral convergence in the number of collocation points $N$
when the boundary shape is a unit circle.

\subsection{Test for ellipse\label{subsec:4.3}}

\begin{figure}
\centering
\includegraphics[width=8.4cm]{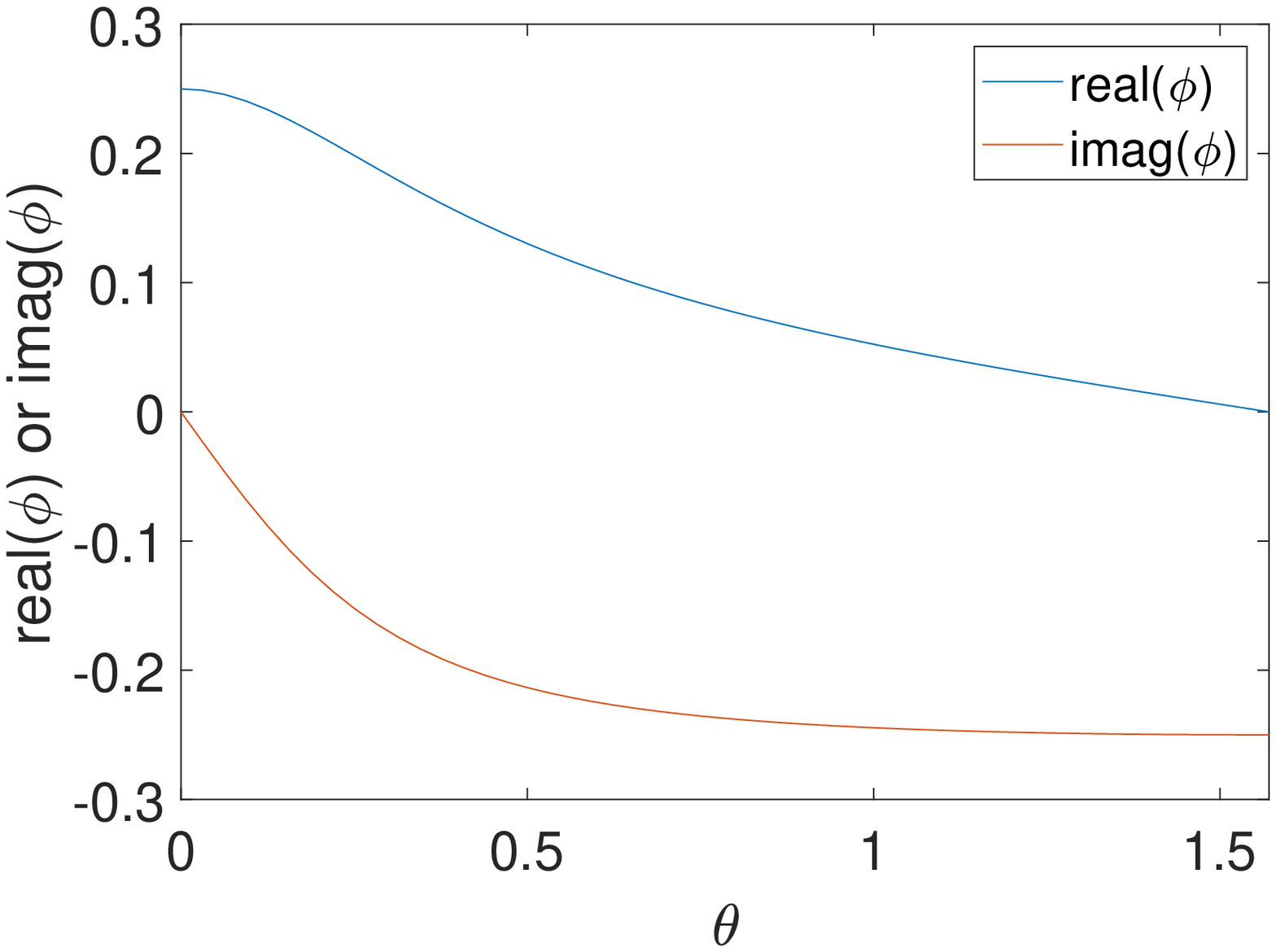}

\caption{\label{fig:5}Real part and imaginary part of exact solution $\varphi_{\mbox{exact}}$
versus $\theta\in[0,\pi/2]$ when $m=0.5$.}
\end{figure}

An ellipse in the complex $z$-plane can be conformally mapped to
the interior circle problem in the $\zeta$-plane with $z(\zeta)=1/\zeta+m\zeta$
for $0<m<1$. Thus, the exact solution for $\chi=0$ is $\varphi(z)=(1-m)/\zeta$
\cite{Muskhelishvili}. The exact solution for the ellipse when $m=0.5$
is shown in Fig.\,\ref{fig:5}.

\begin{figure}
\centering
\includegraphics[width=8.4cm]{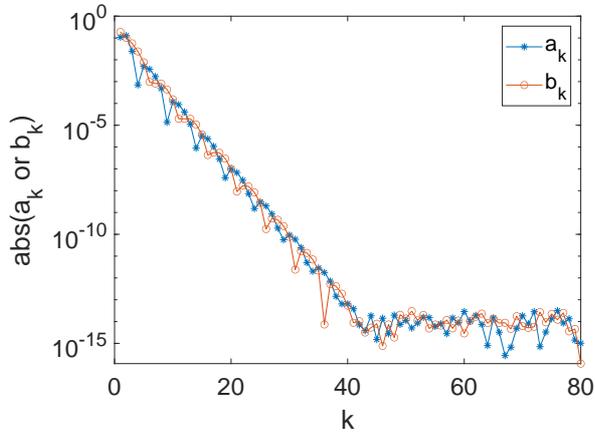}

\caption{\label{fig:6}Result for test case of an ellipse: Absolute value of
coefficients of $a_{k}$ or $b_{k}$ versus $k$ when number of collocation
points $N=80$.}
\end{figure}

\begin{figure}
\centering
\includegraphics[width=8.4cm]{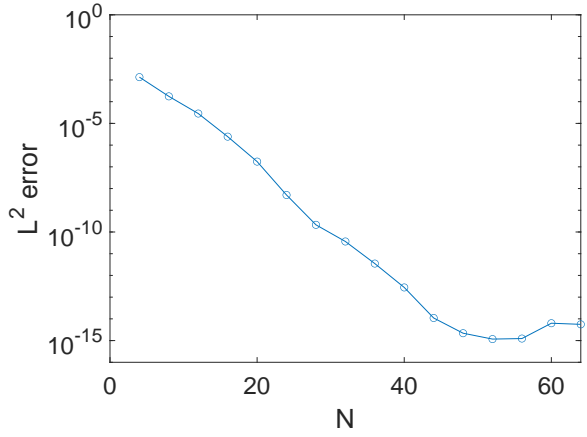}

\caption{\label{fig:7}Result for test case of an ellipse: $L^{2}$ error of
$\varphi(z)$ versus number of collocation points $N$. Here $L^{2}$
error is defined in Section \ref{subsec:4.1}.}
\end{figure}

We test our numerical method for the ellipse with the eccentricity
\begin{equation}
e=\sqrt{1-(1-m)^{2}/(1+m)^{2}}.
\end{equation}
The equation of the ellipse in polar coordinates is
\begin{equation}
r(\theta)=(1-m)/\sqrt{1-(e\cos\theta)^{2}}.
\end{equation}
 We can get $r'(\theta)$ and $r''(\theta)$ by taking the derivatives:
\begin{equation}
r'(\theta)=-\frac{1-m}{2}(1-e^{2}\cos^{2}\theta)^{-3/2}e^{2}\sin2\theta,
\end{equation}
\begin{multline}
r''(\theta)=-\frac{1-m}{2}[-(3/2)\cdot(1-e^{2}\cos^{2}\theta)^{-5/2}\cdot(e^{2}\sin2\theta)^{2}\\
+2(1-e^{2}\cos^{2}\theta)^{-3/2}e^{2}\cos2\theta].
\end{multline}
The calculated Chebyshev coefficients when $m=0.5$ are given in Fig.\,\ref{fig:6}.
Fig.\,\ref{fig:7} shows that the numerical solution converges rapidly
and has error comparable to the truncation error when $N>50$. Thus,
our numerical method converges rapidly for the smooth-boundary test
cases of a circle and an ellipse. 

\subsection{Test for overlapping circles ($\alpha=\pi/3$)\label{subsec:4.4}}

\begin{figure}
\centering
\includegraphics[width=8.4cm]{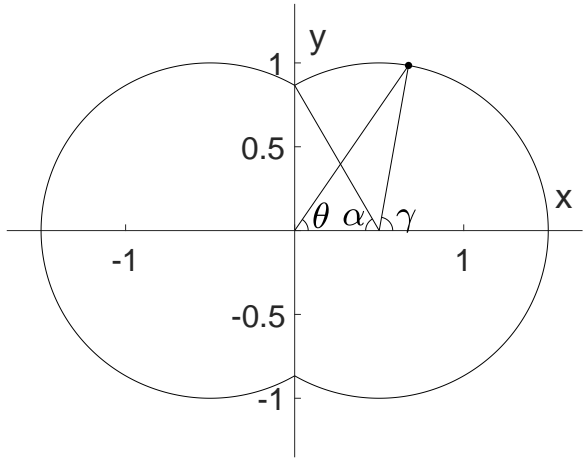}

\caption{\label{fig:8}Overlapping circles shape when $\alpha=\pi/3$.}
\end{figure}

In this test, we show the convergence of the numerical method when
the shape has corners. Consider the shape formed by two overlapping
circles with the same radius (see Fig.\,\ref{fig:8}). The void shape
is given by taking the $x>0$ portion of the unit circle with center
at $(\cos(\alpha),0)$ and reflecting it across the $y$ axis. Thus
the parameter $\alpha$ controls the amount of overlap between the
circles. There are three cases: (i) the degenerate case corresponding
to a single unit circle occurs for $\alpha=\pi/2$; (ii) the \textquotedbl separating
circles\textquotedbl{} case for $0\le\alpha<\pi/2$ as illustrated
in Fig.\,\ref{fig:8}; and (iii) the \textquotedbl collapsing circles\textquotedbl{}
case for $\pi/2<\alpha<\pi$ as illustrated in Fig.\,\ref{fig:12}.
In this test, we consider the \textquotedbl separating circles\textquotedbl{}
case with $\alpha=\pi/3$ for which the void shape has two inward-pointing
corners on the $y$ axis. Since the shape has corners, we apply Eq.\,\eqref{eq:32}
to represent $\varphi(\theta)$ in our numerical method. The order
of the corner term has $\lambda=2.0465$ which is the solution to Eq.\,\eqref{eq:29}
with corner angle $\beta=2\alpha$. Consistent with \cite{williams1952stress},
if the corner angle of the solid is less than $\pi$, there will be
no singularity in the stress at the corner. 

\begin{figure}
\centering
\includegraphics[width=8.4cm]{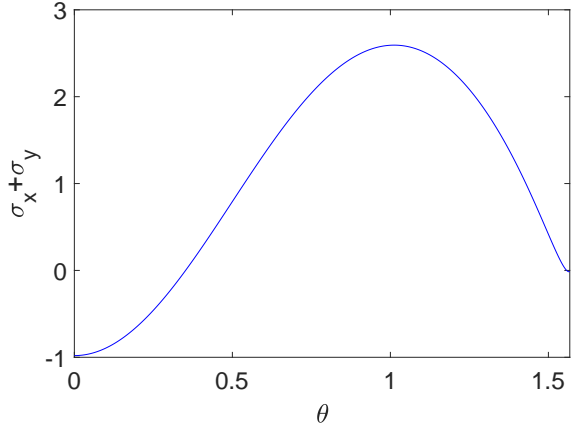}

\caption{\label{fig:9}Exact solution of overlapping circles shape: $\sigma_{x}+\sigma_{y}$
versus $\theta\in[0,\pi/2)$ when $\alpha=\pi/3$ and the far-field
stress parameter $\chi=0$.}
\end{figure}

The exact solution of{\footnotesize{} }$\sigma_{x}+\sigma_{y}$ in
this overlapping circles shape is given in \cite{ling1948stresses}. 
The equation of the shape in polar
coordinates is $r(\theta)=\cos\alpha\cos\theta+\sqrt{1-\sin^{2}\theta\cos^{2}\alpha}$
in the first quadrant. The exact solution for $\sigma_{x}+\sigma_{y}$
is given in \cite{ling1948stresses} by an integral
\begin{multline}
\sigma_{x}+\sigma_{y}=4(\cosh\xi-\cos\alpha)\sin\alpha \, \times\\
\int_{0}^{\infty}\frac{2K-(N_{1}-N_{2})s(s-\cot\alpha\coth s\alpha)}{\sinh2s\alpha+s\sin2\alpha}\times\\
\sinh s\alpha\cos s\xi\,ds,\label{eq:58}
\end{multline}
where $K$ is the solution of the equation
\begin{multline}
4K\int_{0}^{\infty}\frac{\sinh^{2}s\alpha-s^{2}\sin^{2}\alpha}{s(s^{2}+1)(\sinh2s\alpha+s\sin2\alpha)}\,ds\\
+2(N_{1}-N_{2})\int_{0}^{\infty}\frac{s\sin^{2}\alpha}{\sinh2s\alpha+s\sin2\alpha}\,du=N_{1},
\end{multline}
where $\xi$ is defined by
\begin{equation}
\cosh\xi=\frac{1+\cos\alpha\cos\gamma}{\cos\alpha+\cos\gamma},
\end{equation}
$\gamma=\theta+\arcsin(\sin\theta\cos\alpha)$ is the center angle
(see Fig.\,\ref{fig:8}), $N_{1}$ is the tension parallel to the
$x$-axis and $N_{2}$ is the tension parallel to the $y$-axis. The value
of $\sigma_{x}+\sigma_{y}$ for the exact solution is evaluated using the
MATLAB numerical integration function 'integral' \cite{shampine2008vectorized}.
The exact solution for $\sigma_x + \sigma_y$ is shown in Fig.\,\ref{fig:9} for
the overlapping circles case with $\alpha=\pi/3$ and longitudinal tension ($N_{1}=1$
and $N_{2}=0$, corresponding to $\chi = 0$). 
The accuracy of the numerical integration near the corner has been verified
by comparing to the asymptotic behavior of the integral near the corner.  
Derivation of the asymptotic behavior of the integral is given in 
Appendix A.

\begin{figure}
\centering
\includegraphics[width=8.4cm]{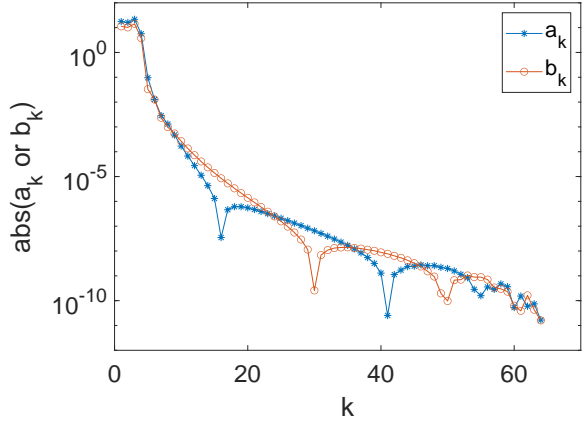}

\caption{\label{fig:10}Overlapping circles with $\alpha=\pi/3$ and $\chi=0$:
Absolute value of coefficients of $a_{k}$ or $b_{k}$ versus $k$
when number of collocation points $N=64$. }
\end{figure}
\begin{figure}
\centering
\includegraphics[width=8.4cm]{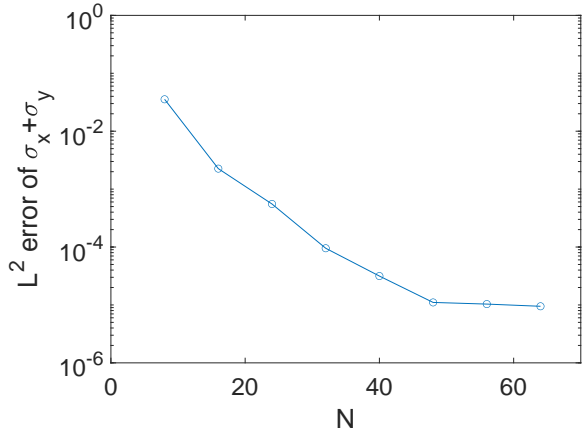}

\caption{\label{fig:11}Overlapping circles with $\alpha=\pi/3$ and $\chi=0$:
$L^{2}$ error of $\sigma_{x}+\sigma_{y}$ versus number of collocation
points $N$.}
\end{figure}

We test our numerical method for the overlapping circles case $\alpha = \pi/3$
for the uniaxial tension case $\chi = 0$.
Fig.\,\ref{fig:10} and Fig.\,\ref{fig:11} demonstrate that our
numerical method works when the shape of the hole has a corner
with angle greater than $\pi$. Fig.\,\ref{fig:10} shows that 
the magnitude of the coefficients
$a_{k}$ and $b_{k}$ decays to $10^{-10}$ for $N=64$ and Fig.\,\ref{fig:11} 
shows that the $L^{2}$ error
of $\sigma_{x}+\sigma_{y}$ is about $10^{-5}$ when $N$ is
larger than about 45.  While the result is not as accurate as
the results for smooth shapes in Section \ref{subsec:4.2} and Section
\ref{subsec:4.3}, the accuracy is still very good for relatively small $N$.
See Section \ref{subsec:4.7} for a discussion
of the error.

\subsection{Test for overlapping circles ($\alpha=2\pi/3$)}

\begin{figure}
\centering
\includegraphics[width=8.4cm]{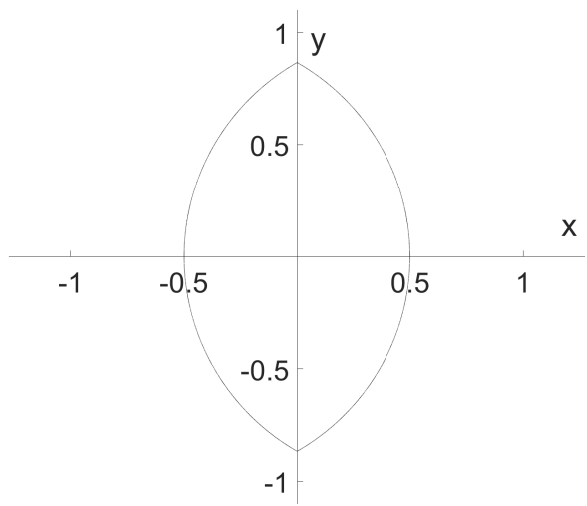}

\caption{\label{fig:12}Overlapping circles shape when $\alpha=2\pi/3$.}
\end{figure}
\begin{figure}
\centering
\includegraphics[width=8.4cm]{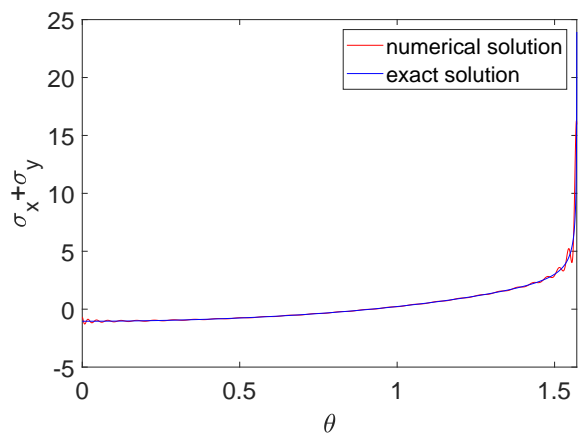}

\caption{\label{fig:13}Overlapping circles case ($\alpha=2\pi/3$), no corner
term:{\footnotesize{} }Numerical solution and exact solution of $\sigma_{x}+\sigma_{y}$
versus $\theta\in[0,\pi/2)$. Parameters are $N=64$, $\chi=0$.}
\end{figure}
\begin{figure}
\centering
\includegraphics[width=8.4cm]{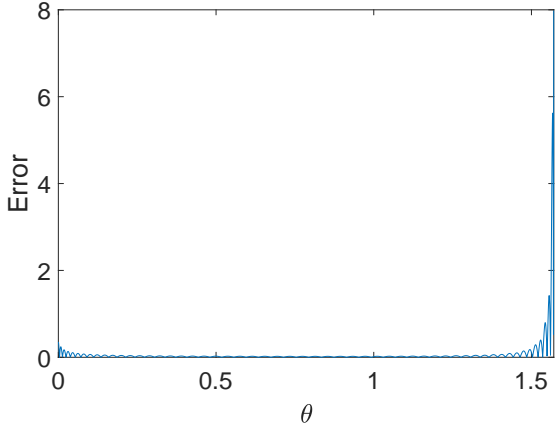}

\caption{\label{fig:14}Overlapping circles case ($\alpha=2\pi/3$), no corner
term: Error of $\sigma_{x}+\sigma_{y}$ versus $\theta\in[0,\pi/2)$.
Parameters are $N=64$, $\chi=0$.}

\end{figure}

We show the effect of the corner term in Eq.\,\eqref{eq:32} in this
test. The shape of "collapsing'' overlapping circles when $\alpha=2\pi/3$ is shown
in Fig.\,\ref{fig:12}. When the hole shape has corners with corner
angle less than $\pi$, the stress has a singularity at the corner
\cite{williams1952stress}. If we use the expansion \eqref{eq:31}
for $\varphi(\theta)$ without the corner term the result is as shown
in Figs.\,\ref{fig:13} and \ref{fig:14}. The error of $\sigma_{x}+\sigma_{y}$
is large near the corner at $\theta = \pi/2$. Expansion \eqref{eq:31} uses
a set of polynomials to approximate a $\varphi(\theta)$ that has a discontinuous
derivative at the end of the interval. Thus, oscillations are expected
near the point of discontinuity (Runge's phenomenon \cite{epperson1987runge}).


\begin{figure}
\centering
\includegraphics[width=8.4cm]{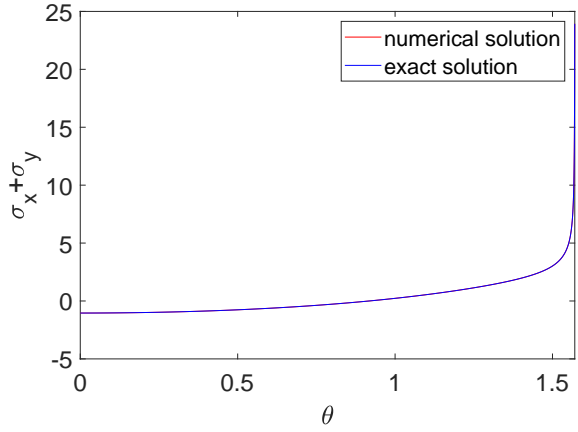}

\caption{\label{fig:15}Overlapping circles case ($\alpha=2\pi/3$) with corner
terms:{\footnotesize{} }Numerical solution and exact solution of $\sigma_{x}+\sigma_{y}$
versus $\theta\in[0,\pi/2)$. Parameters are $N=64$, $\chi=0$.}
\end{figure}

\begin{figure}
\centering
\includegraphics[width=8.4cm]{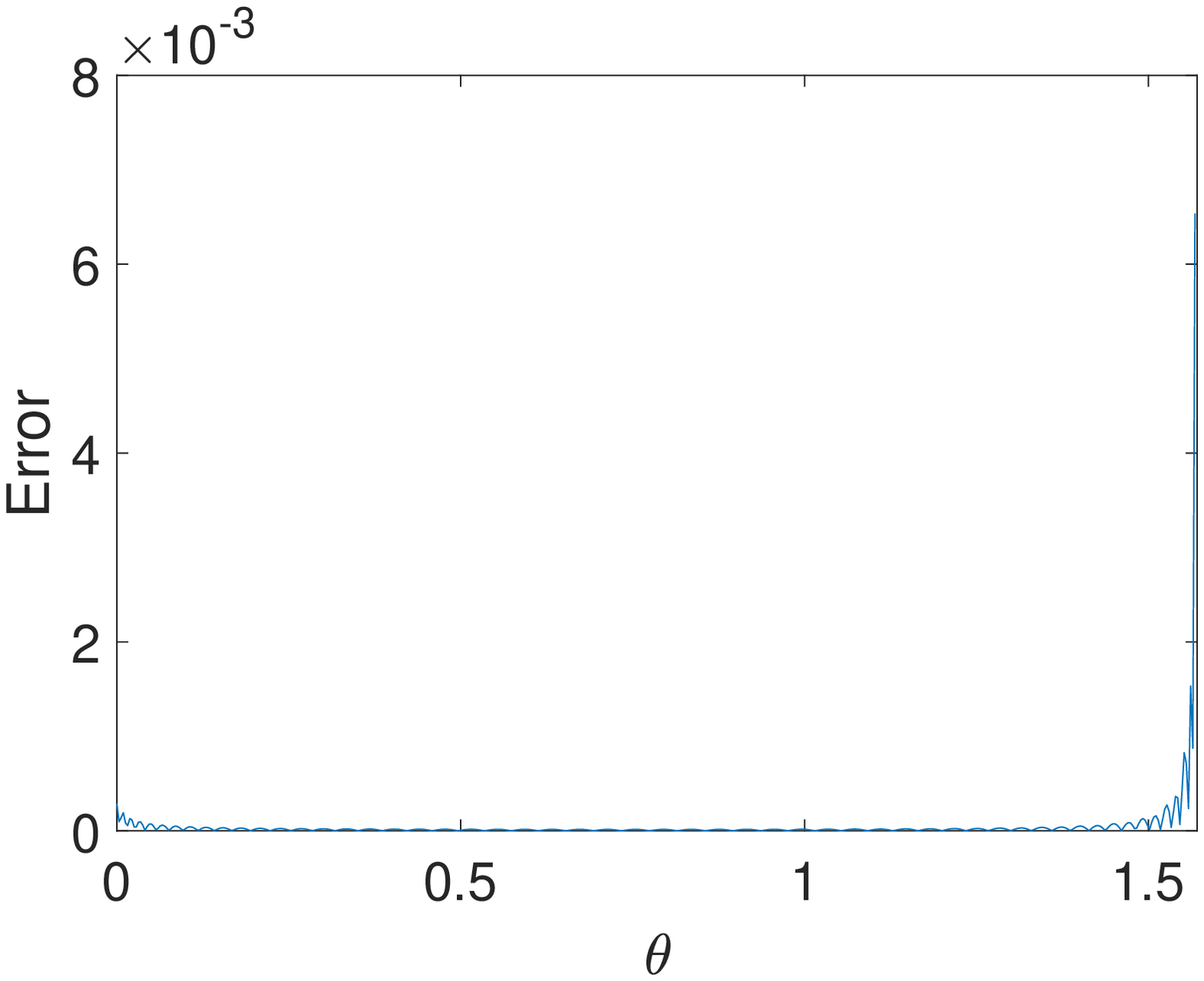}

\caption{\label{fig:16}Overlapping circles case ($\alpha=2\pi/3$) with corner
terms: Error of $\sigma_{x}+\sigma_{y}$ versus $\theta\in[0,\pi/2)$.
Parameters are $N=64$, $\chi=0$.}
\end{figure}

\begin{figure}
\centering
\includegraphics[width=8.4cm]{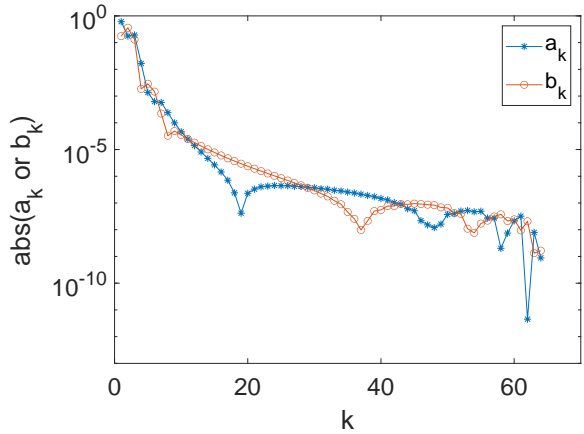}

\caption{\label{fig:17}Overlapping circles case ($\alpha=2\pi/3$) with corner
terms: Absolute value of coefficient of $a_{k}$ or $b_{k}$ versus
index $k$. 
Parameters are $N=64$, $\chi=0$.}
\end{figure}

\begin{figure}
\centering
\includegraphics[width=8.4cm]{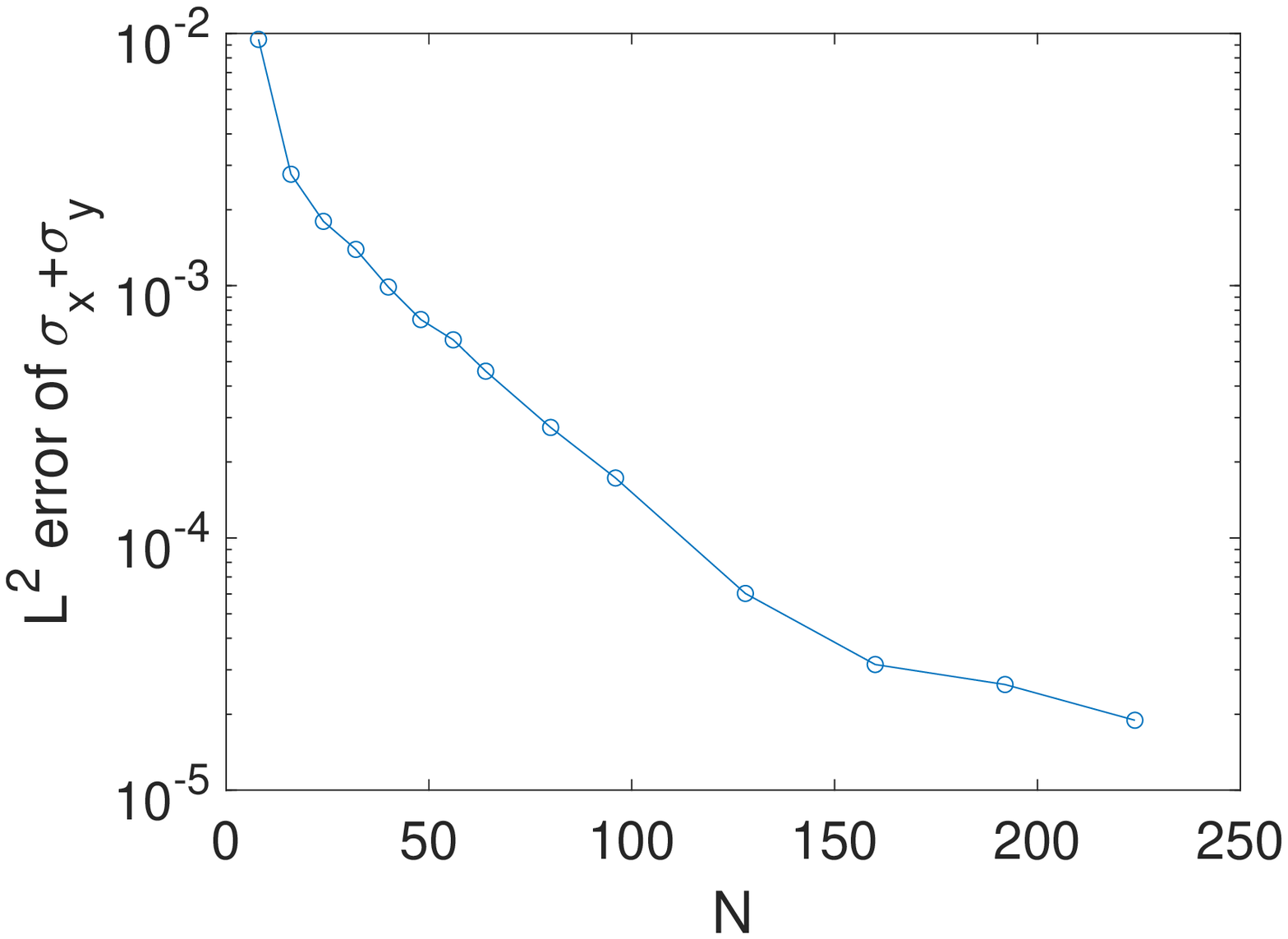}

\caption{\label{fig:18}Overlapping circle case ($\alpha=2\pi/3$) with corner
terms: $L^{2}$ error of $\sigma_{x}+\sigma_{y}$ versus number of
collocation points $N$ for $\chi=0$.}
\end{figure}

We can reduce the error by introducing the corner term in $\varphi(\theta)$
as given in Eq.\,\eqref{eq:32}. The order of the corner term $\lambda-1$
can be obtained by Eq.\,\eqref{eq:29} using the known value of the
corner angle $\beta$. Example numerical results for $N=64$ are shown in Figs.\,\ref{fig:15}-\ref{fig:17}.
Comparing the magnitude of the error to the magnitude of the solution near the
corner we note that the maximum of the relative error is less than $10^{-3}$ approaching the corner.
Fig.\,\ref{fig:18}
shows the $L^{2}$ error of $\sigma_{x}+\sigma_{y}$ versus number
of collocation points.  The $L^{2}$ error is less than $10^{-4}$ when $N$ is
large, which is comparable to the error magnitude for the nonsingular
case.
So, even in this case with a stress singularity
near the corner, the numerical solution is of good accuracy for moderately large $N$.


\subsection{Contour plots of stresses}

\begin{figure}
\centering
\includegraphics[width=8.4cm,height=6.3cm]{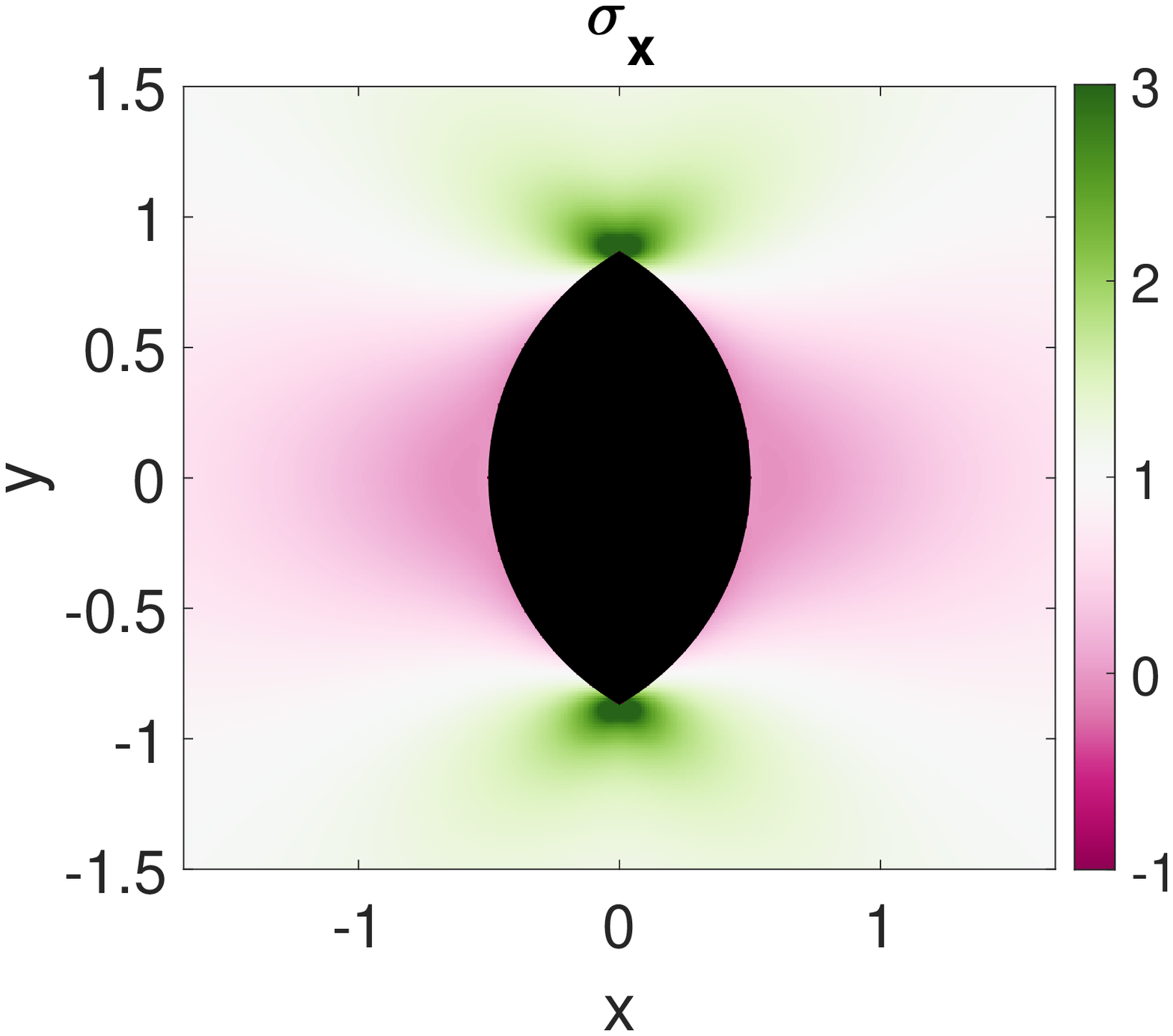}

\includegraphics[width=8.4cm,height=6.3cm]{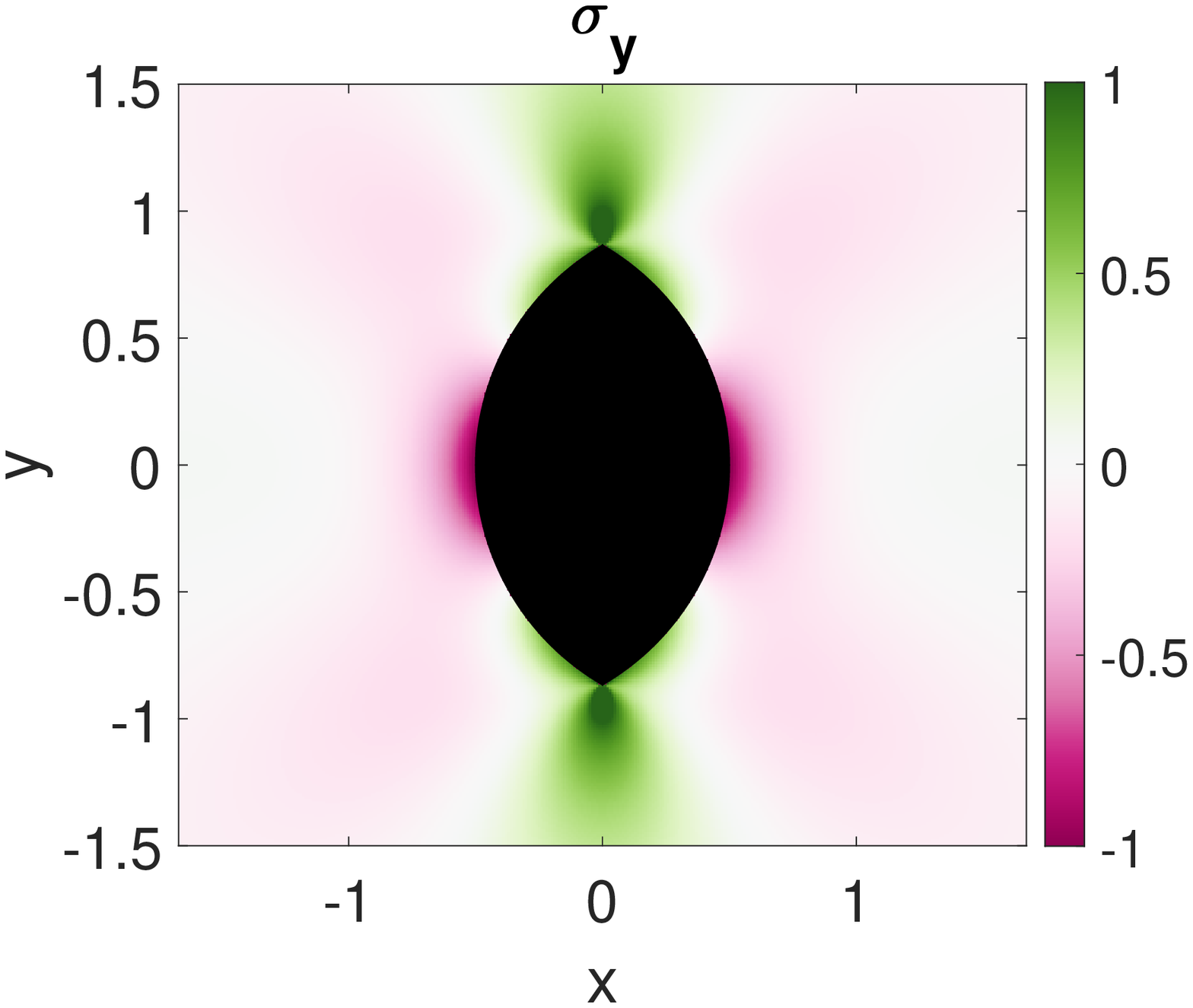}

\includegraphics[width=8.4cm,height=6.3cm]{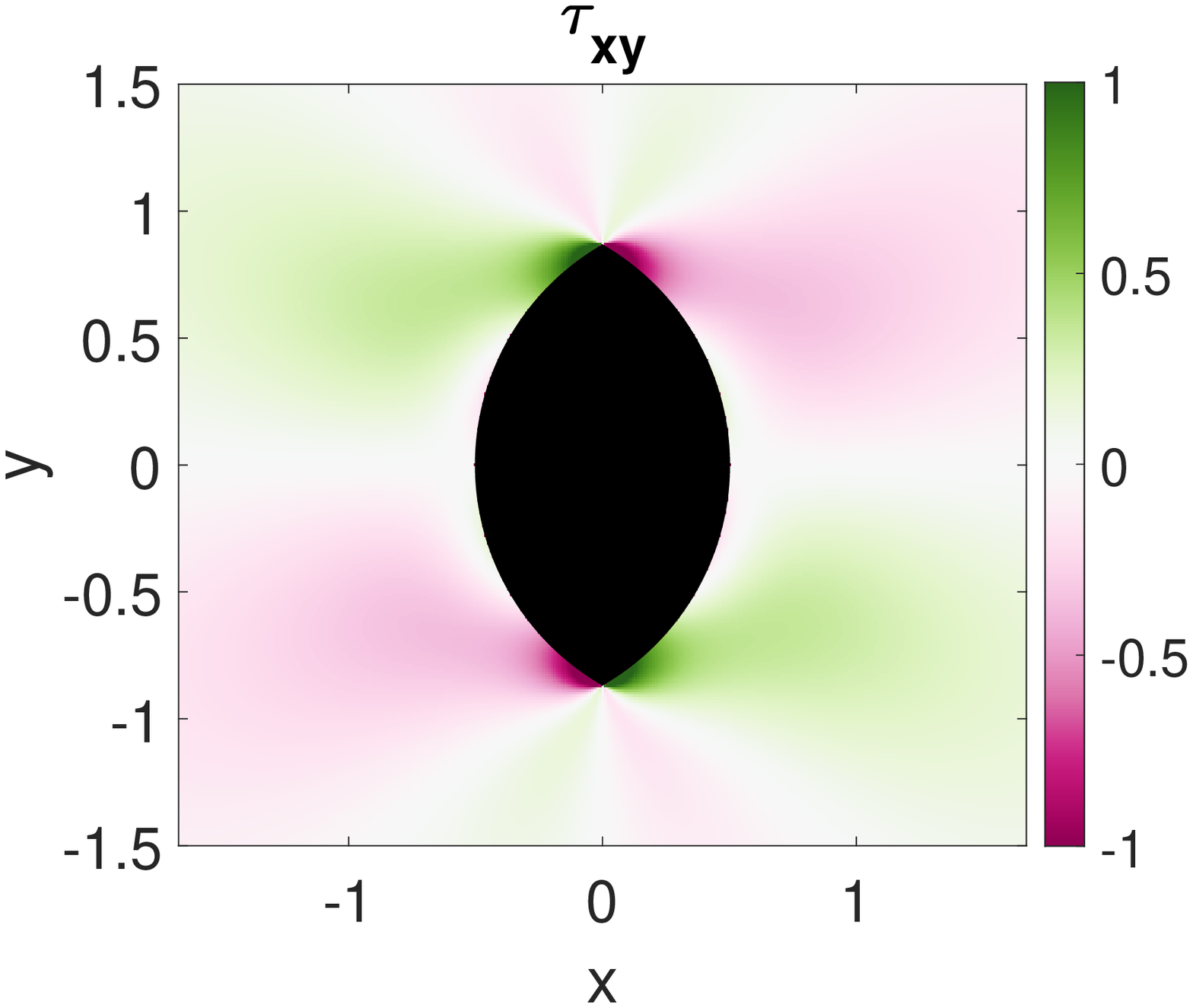}

\caption{\label{fig:19}Distribution of stresses for overlapping circles case
with $\alpha=2\pi/3$, $\chi=0$. Note
the exact stress may exceed the range of the color bar due to the
singularity at the corner.}
\end{figure}

\noindent Our numerical method gives the boundary value of $\varphi$
with small error. We can find the boundary value of the analytic function
$h(z)$ on $\partial D$ with small error by Eq.\,\eqref{eq:12}.
Then analytic functions $\varphi,h$ can be extended to any $\zeta$
in $D$ by Cauchy's integral formula:
\begin{equation}
\varphi(\zeta)=\frac{1}{2\pi i}\underset{-L}{\int}\frac{\varphi(z)}{z-\zeta}\,dz,\quad h(\zeta)=\frac{1}{2\pi i}\underset{-L}{\int}\frac{h(z)}{z-\zeta}\,dz.
\end{equation}
The integral is evaluated using by MATLAB numerical integration function since there is no 
singularity if $\zeta$ is not on the boundary. 
From $\varphi(\zeta)$ and $h(\zeta)$ in $D$ we can determine the stresses at any point in $D$ using Eqs.\,\eqref{eq:5}-\eqref{eq:6}, where we
apply the finite difference method with small 
step size to evaluate the derivatives of $\varphi$ and $h$.
Fig.\,\ref{fig:19} illustrates the results
for the stress distribution for the overlapping circles case. 

\subsection{Error analysis\label{subsec:4.7}}

\noindent From our numerical results for the cases of the circle and
ellipse, we see in these smooth boundary cases our numerical method
gives spectral convergence. In this case, both real and imaginary
parts of the Goursat function are smooth and thus can be well approximated
by Chebyshev series. Nested Gaussian quadrature achieves spectral
convergence for the numerical integration, and analyticity constraints
Eq.\,\eqref{eq:51} guarantee that $\varphi$ must correspond to
boundary values of an analytic function which ensure the correctness
and uniqueness of the solution. 

In the case with the non-smooth shape, the numerical solution is close
to the exact solution with $L^2$ error less than $10^{-4}$ even
when the stresses have singularities near the corner. 
We argue that this error is a combination of not completely resolving the
non-analytic behavior near the corner and poor conditioning of the discretized
system.  
The nonanalytic term in the expansion Eq.\,\eqref{eq:32} was derived from
asymptotic analysis near the corner. While we capture the the dominant
behavior of the singularity with the corner term in Eq.\,\eqref{eq:32},
we do not capture other non-integer powers in the expansion of
the wedge solution \cite{williams1952stress}, and thus
we would not expect to recover the exact solution
near the corner.  
These other non-integer powers could be determined from the higher-order
solutions to Eq.\,\eqref{eq:29} (which is the same as Eq.\,\eqref{eq:76} 
that could be
used to find the sub-dominant terms in
our asymptotic expansion of the exact solution in Appendix A).

Another contribution to the error is due to poor conditioning.  
Since the corner term in Eq.\,\eqref{eq:32} is
not orthogonal to Chebyshev polynomials, there is some redundancy
in the expansion which makes the matrix problem from collocation method
somewhat ill-conditioned due to the lack of independence of the corner
coefficient $a_{N-1}$ and the coefficients of the regular Chebyshev
expansion. It may be possible in a future work to use the wedge solution
to give not only the power but also the coefficient of this corner
term from the asymptotic analysis. In this case, by determining the
$a_{N-1}$ coefficient from the asymptotic analysis, the resulting regular
Chebyshev expansion would be used to fit $\varphi-a_{N-1}(\pi/2-\theta)^{\lambda-1}$
which would be better conditioned. If we have to retain the corner
term as a term with an unknown coefficient, then to accurately reproduce
the non-integer power series in the vicinity of the corner, and remove
the ill-conditioning of the current matrix, we could consider developing
an orthogonal basis based on the non-integer power from the corner
singularity. 

\section{Summary}

We developed a numerical method to determine the elastic stresses around
a hole in an infinite plate. Our numerical method is based on the
boundary integro-differential equations obtained using complex Goursat
functions. We represent the real and imaginary part of the boundary value
of the complex Goursat function $\varphi$ by a series of Chebyshev
polynomials with possible corner terms. The corner term is chosen
based on the asymptotic analysis of the stress near the corner and
captures the dominant (possibly singular) behavior of the stress.
The boundary integro-differential equation is solved numerically using the
collocation method at Gauss-Legendre collocation points. To improve
the accuracy of the numerical method, we separate the singularity
in the integro-differential equation and evaluate the remaining parts
by nested Gaussian quadrature. The boundary equation for $\varphi$
is augmented by analyticity constraints to ensure
$\varphi$ is an analytic function in the solid. Then we solve the
over-determined system of the integro-differential equations, end
of the interval conditions together with analyticity constraints to
determine the boundary value of $\varphi$. For the case of a smooth
boundary shape, our numerical method converges to the exact solution
spectrally with a small number of collocation points. For the case
when the boundary shape has a corner, our numerical method includes
a corner term derived from asymptotic analysis of the stress singularity
at the corner. In this case we obtain elastic stresses accurate to
relative error less than $10^{-3}$. Finally, we obtain the 
distribution of stresses 
in the solid by taking the Cauchy integral of the Goursat
function on the boundary of the hole. 

In the present work, we make use of the corner angle to determine
the asymptotic behavior of the stress singularity in the vicinity
of a corner and use this asymptotic behavior to modify our numerical
method. While the results presented here demonstrate the approach
for only a single corner in our computational domain, the ideas could
be extended to handle multiple corners in a straightforward way. Also,
we make use of an assumed two-fold symmetry for the shape to reduce
the number of unknowns for the integral equation for the Goursat function
$\varphi$ and to eliminate the homogenous solution. For holes with
asymmetric shape these symmetry assumptions could be relaxed. In the
asymmetric case the shape and the Goursat function could be represented
on the interval $0\le\theta\le2\pi$ using an appropriate spectral
representation.

The independence of analyticity constraints in Section \ref{subsec:3.4}
deserves further attention. Eq.\,\eqref{eq:32} at collocation points
adds $2N-2$ equations to the linear system to ensure the analyticity
which causes the linear system to be over-determined. If we formulate
the integral equations in a way that the solution for $\varphi$ is
guaranteed to be the boundary value of an analytic function (see,
for example \cite{mikhlin1957integral}), the size of the linear system
can be reduced to about half as large as before which would speed up the numerical
method considerably. 

\section*{Appendix A: Asymtotic analysis near the corner for overlapping circles
case.}

Here we derive the asymptotic behavior of the stress near the corner for the 
overlapping circle case of Section 4.5. The trace of the stress tensor for 
the overlapping circles case is given by Eq.\,\eqref{eq:58}
which we write as
\begin{equation}
\sigma_{x}+\sigma_{y}=H(\xi,\alpha)\int_{0}^{\infty}F(s,\alpha)\cos s\xi\,ds,\label{eq:62}
\end{equation}
where 
\begin{equation}
H(\xi,\alpha)=4(\cosh\xi-\cos\alpha)\sin\alpha,
\end{equation}
\begin{equation}
F(s,\alpha)=\frac{2K-(N_{1}-N_{2})s(s-\cot\alpha\coth s\alpha)}{\sinh2s\alpha+s\sin2\alpha}\cdot\sinh s\alpha
\end{equation}
with 
\begin{equation}
\cosh\xi=\frac{1+\cos\alpha\cos\gamma}{\cos\alpha+\cos\gamma}
\end{equation}
and
\begin{equation}
\gamma=\theta+\arcsin(\sin\theta\cos\alpha).\label{eq:66}
\end{equation}
Here, the polar angle from the center of the circle $\gamma$, polar angle
$\theta$, amount of overlap between the circles $\alpha$, tensions
$N_{1}$ and $N_{2}$ are defined as in Section \ref{subsec:4.4}. We use
$\alpha=2\pi/3$ as an example in this Appendix. The result can be generalized to other 
$\alpha$ in $\pi/2<\alpha<\pi$ using a similar approach as described here.
The corner location is at $\theta=\pi/2$. As $\theta\rightarrow\pi/2$,
$\gamma\rightarrow\pi/2+\arcsin(\cos\alpha)$. Thus, 
\begin{equation}
\cos\gamma\rightarrow\cos\left[\pi/2+\arcsin(\cos\alpha)\right]=-\cos\alpha.
\end{equation}
So 
\begin{equation}
\cosh\xi=\frac{1+\cos\alpha\cos\gamma}{\cos\alpha+\cos\gamma}\rightarrow\infty,\label{eq:68}
\end{equation}
which gives $\xi\rightarrow\infty$ as $\theta\rightarrow\pi/2$.
Thus the $\cos(s\xi)$ term 
in integral \eqref{eq:62} is highly oscillatory
near the corner which may cause inaccuracy of the numerical integration. 

We determine the behavior of the integral near the corner by applying
asymptotic analysis. Let $\theta=\pi/2-\varepsilon$ where $0<\varepsilon\ll1$
then $\cos\theta=\sin\varepsilon=\varepsilon-\varepsilon^{3}/6+\mathcal{O}(\varepsilon^{5})$
and $\sin\theta=\cos\varepsilon=1-\varepsilon^{2}/2+\mathcal{O}(\varepsilon^{2})$.
By Eq.\,\eqref{eq:66}, we have
\begin{equation}
\cos\gamma=\cos\theta\sqrt{1-\sin^{2}\theta\cos^{2}\alpha}-\sin^{2}\theta\cos\alpha.
\end{equation}
Thus we have
\begin{equation}
\cos\gamma=-\cos\alpha+\varepsilon\sin\alpha+\mathcal{O}(\varepsilon^{2}).\label{eq:70}
\end{equation}
Substitute Eq.\,\eqref{eq:70} into Eq.\,\eqref{eq:68}, the asymptotic
approximation of $\cosh\xi$ is
\begin{equation}
\cosh\xi=\frac{e^{\xi}+e^{-\xi}}{2}=\frac{1}{\varepsilon}\sin\alpha+\cos\alpha+\mathcal{O}(\varepsilon).\label{eq:71}
\end{equation}
Thus, 
\begin{equation}
\xi=\ln\left[\frac{2}{\varepsilon}\sin\alpha+2\cos\alpha+\mathcal{O}(\varepsilon)\right], \label{eq:72}
\end{equation}
Now consider the integral in \eqref{eq:62}. $F$ has properties
$F(-s)=F(s)$, $F(s)\rightarrow0$ as $s\rightarrow\pm\infty$ and
$F(s,\alpha)$ is bounded for all $s$ ($s=0$ is a removable singularity),
so we can rewrite the integral on $-\infty<s<+\infty$ as
\begin{equation}
I=\int_{0}^{\infty}F(s,\alpha)\cos s\xi\,ds=\frac{1}{2}\int_{-\infty}^{\infty}F(s,\alpha)\cos s\xi\,ds.
\end{equation}
We evaluate $I$ by considering the contour integral on complex plane
$z=s+it$:
\begin{equation}
\tilde{I}=\underset{R\rightarrow\infty}{\lim}\int_{C}F(z,\alpha)e^{i\xi z}\thinspace dz, \label{eq:74}
\end{equation}
\begin{equation}
I=\frac{1}{2}\mbox{Re}\{\tilde{I}\},
\end{equation}
where $C$ is the line along $s$ axis from $-R$ to $-R$.
$F(z,\alpha)$ is analytic in upper half plane $\mbox{Im}(z)>0$ except
zeros of the denominator located by the roots of
\begin{equation}
\sinh2z\alpha+z\sin2\alpha=0.\label{eq:75}
\end{equation}
\begin{figure}
\centering
\includegraphics[width=8.4cm]{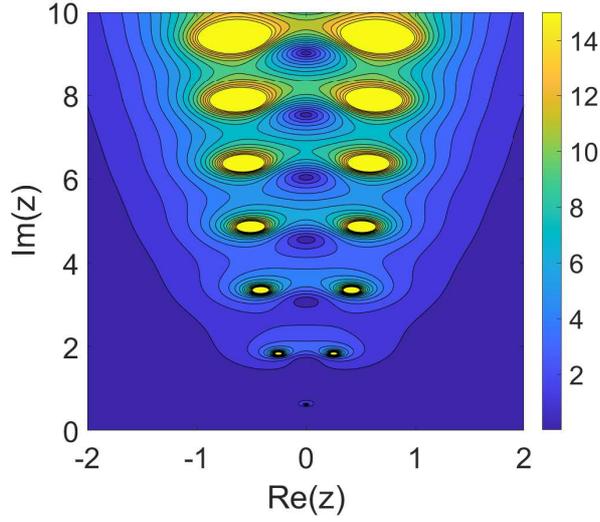}
\caption{\label{fig:20}Contour plot of $\left|F\right|$ on upper half plane for $\alpha=2\pi/3$.}
\end{figure}
There are infinitely many singularities of $F(z,\alpha)$ on the upper half plane (see Fig.\,\ref{fig:20}). 
\begin{figure}
\centering
\includegraphics[width=8.4cm]{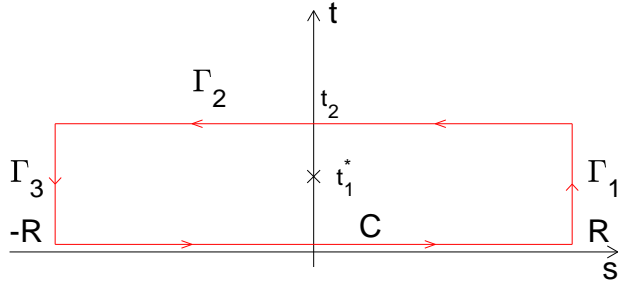}
\caption{\label{fig:21}Complex contour integral.}
\end{figure}
We claim that since $\tilde{I}$ contains $e^{i\xi z}=e^{-\xi t}e^{i\xi s}$ and $\xi\gg 1$, the integral $\tilde{I}$ 
can be approximated by the contribution from the residual of the singularity that occurs at the location with smallest imaginary 
part $t$ in the upper half plane. For $\alpha=2\pi/3$, as shown in Fig.\,\ref{fig:20}, this first singularity lies on the
imaginary axis, and so satisfies $s=0$ and 
\begin{equation}
\sin2\alpha t+t\sin2\alpha=0.\label{eq:76}
\end{equation}
For $\alpha=2\pi/3$ the solution of Eq.\,\eqref{eq:76} is $t_1^*\approx 0.6157$ and the first singularity is at $z=it_1^*$. 
We construct a rectangular contour on the complex plane with the first singularity inside as in Fig.\,\ref{fig:21} to evaluate 
the integral \eqref{eq:74}, where $C$ is the line segment from $z=-R$ to $z=R$, 
$\Gamma_1$ is the line segment from $z=R$ to $z=R+it_2$, $\Gamma_2$ is the line segment from $z=R+it_2$ to $z=-R+it_2$ 
and $\Gamma_3$ is the line segment from $z=-R+it_2$ to $z=-R$. The height of the rectangle is defined by $t_2=\pi/(2\alpha)>t_1^*$ such that only 
the first singularity is inside the contour. By the residue theorem and taking the limit $R\rightarrow\infty$, 
\begin{equation}
\underset{R\rightarrow\infty}{\lim}\int_{C+\Gamma_1+\Gamma_2+\Gamma_3}F(z,\alpha)e^{i\xi z}\thinspace dz
=2\pi i\cdot\mbox{ residue of }F(z,\alpha)e^{i\xi z}\mbox{ at }z=it_{1}^{*}.
\end{equation}
Now consider the integrals on $\Gamma_1$, $\Gamma_2$ and $\Gamma_3$. First, on $\Gamma_1$
\begin{multline}
\underset{R\rightarrow\infty}{\lim}\left|\int_{\Gamma_1}F(z,\alpha)e^{i\xi z}\thinspace dz\right|=\underset{R\rightarrow\infty}{\lim}\left|\int_0^{t_2}F(R+it,\alpha)e^{-\xi t+i\xi R}\thinspace dt\right|\\
\leq\int_0^{t_2}\underset{R\rightarrow\infty}{\lim}\left|F(R+it,\alpha)e^{-\xi t}\right|\thinspace dt\\
\leq\int_0^{t_2}\underset{R\rightarrow\infty}{\lim}\Big[\frac{(2\left|K\right|+\left|N_1-N_2\right|\left|z\right|^2)\left|\sinh(\alpha z)\right|}{|\sinh(2\alpha z)|-|z|\left|\sin{2\alpha}\right|}
+\frac{\left|N_1-N_2\right|\left|\cot\alpha\right|\left|z\right|\left|\cosh(\alpha z)\right|}{|\sinh(2\alpha z)|-|z|\left|\sin{2\alpha}\right|}\Big]\thinspace dt.
\end{multline}
We can bound the $|F|$ using following inequalities
\begin{equation}
|\sinh\alpha z|=|\sinh\alpha R\cos\alpha t+i\cosh\alpha R\sin\alpha t|
\leq 2|\cosh\alpha R|=|e^{\alpha R}+e^{-\alpha R}|,
\end{equation}
\begin{equation}
|\cosh\alpha z|=|\cosh\alpha R\cos\alpha t+i\sinh\alpha R\sin\alpha t|
\leq 2|\cosh\alpha R|=|e^{\alpha R}+e^{-\alpha R}|,
\end{equation}
\begin{equation}
|z|\leq \sqrt{R^2+t_2^2},
\end{equation}
\begin{multline}
|\sinh 2\alpha z|=|\cosh 2\alpha R\cos 2\alpha t+i\sinh 2\alpha R\sin 2\alpha t|\\
\geq |\sinh 2\alpha R\cos 2\alpha t+i\sinh 2\alpha R\sin 2\alpha t|\geq |\sinh 2\alpha R|.
\end{multline}
Then
\begin{multline}
\underset{R\rightarrow\infty}{\lim}\left|F\right|\leq\underset{R\rightarrow\infty}{\lim}\Big[\frac{(2\left|K\right|+\left|N_1-N_2\right|(R^2+t_2^2)\left|2\cosh(\alpha R)\right|}{|\sinh(2\alpha R)|-\sqrt{R^2+t_2^2}\left|\sin{2\alpha}\right|}\\
+\frac{\left|N_1-N_2\right|\left|\cot\alpha\right|\sqrt{R^2+t_2^2}\left|2\cosh(\alpha R)\right|}{|\sinh(2\alpha R)|-\sqrt{R^2+t_2^2}\left|\sin{2\alpha}\right|}\Big]=0.\\
\end{multline}
Thus, 
\begin{equation}
\underset{R\rightarrow\infty}{\lim}\int_{\Gamma_1}F(z,\alpha)e^{i\xi z}\thinspace dz=0.
\end{equation}
For the same reason, 
\begin{equation}
\underset{R\rightarrow\infty}{\lim}\int_{\Gamma_3}F(z,\alpha)e^{i\xi z}\thinspace dz=0.
\end{equation}
For the $\Gamma_2$ contour
\begin{multline}
\underset{R\rightarrow\infty}{\lim}\left|\int_{\Gamma_2}F(z,\alpha)e^{i\xi z}\thinspace dz\right|
\leq\underset{R\rightarrow\infty}{\lim}\int_{-R}^{R}\left|F(s+i{t_2},\alpha)e^{-\xi t_2+i\xi s}\right|\thinspace ds\\
\leq \Big[2\left|K\right|\left|\int_{-\infty}^{\infty}\frac{\sinh(\alpha(s+it_2))}{\sinh(2\alpha(s+it_2))+(s+it_2)\sin(2\alpha)}\thinspace ds\right|\\
+\Big(\left|\int_{-\infty}^{\infty}\frac{(s+it_2)^2\sinh(\alpha(s+it_2))}{\sinh(2\alpha(s+it_2))+(s+it_2)\sin(2\alpha)}\thinspace ds\right|\\
+\left|\int_{-\infty}^{\infty}\frac{(s+it_2)\cot{\alpha}\cosh(\alpha(s+it_2))}{\sinh(2\alpha(s+it_2))+(s+it_2)\sin(2\alpha)}\thinspace ds \right| \Big) 
\cdot\left|N_1-N_2\right|\Big]\times e^{-\xi t_2}.\label{eq:87}
\end{multline}
All three integrals in the square brackets of Eq.\,\eqref{eq:87} are finite, so the bound on the integral on contour $\Gamma_2$ is of
order $e^{-\xi t_2}$ for $\xi\gg 1$. 

Finally, consider the residue from the singularity:
\begin{multline}
2\pi i\cdot\mbox{ residue of }F(z,\alpha)e^{i\xi z}\mbox{ at }it_{1}^{*}
=2\pi i\underset{z\rightarrow z_{1}^{*}}{\lim}(z-z_{1}^{*})F(z,\alpha)e^{i\xi z}\\
=-2\pi e^{-\xi t_{1}^{*}}\cdotp\frac{2K+(N_{1}-N_{2})t_{1}^{*}(t_{1}^{*}-\cot\alpha\coth\alpha t_{1}^{*})}{2\alpha\cos2\alpha t_{1}^{*}+\sin2\alpha}
\times\sin\alpha t_{1}^{*}\sim \mathcal{O}(e^{-\xi t_1^*}).
\end{multline}
From Eq.\,\eqref{eq:72}, $\xi\gg 1$ near the corner. Then the integral on $\Gamma_2$ is asymptotically smaller than the residue at $z=it_1^*$ because $e^{-\xi t_2}\ll e^{-\xi t_1^*}$. 
Thus the dominant asymptotic contribution to the integral is
\begin{equation}
I\sim\frac{1}{2}\mbox{Re}\{2\pi i\cdot\mbox{ residue of }F(z,\alpha)e^{i\xi z}\mbox{ at }z=it_1^*\}.
\end{equation}
Using this result in Eq.\,\eqref{eq:72} we obtain
\begin{multline}
\sigma_{x}+\sigma_{y}=H(\xi,\alpha)\cdot I
=4(\cosh\xi-\cos\alpha)\sin\alpha\cdot I\\
\sim-\frac{2K+(N_{1}-N_{2})t_{1}^{*}(t_{1}^{*}-\cot\alpha\coth\alpha t_{1}^{*})}{2\alpha\cos2\alpha t_{1}^{*}+\sin2\alpha} 
2\pi\sin\alpha\sin\alpha t_{1}^{*}\cdot(\frac{2}{\varepsilon}\sin\alpha)^{1-t_{1}^{*}},
\end{multline}
where $t_{1}^{*}$ is the smallest nonzero root of Eq.\,\eqref{eq:76}.
Recalling that $\varepsilon=\pi/2-\theta$ is the proximity to the
corner, since  $t_{1}^{*}<1$, $\sigma_{x}+\sigma_{y}$ has an integrable singularity at the corner.
Note the exponent $1-t_1^*$ matches the exponent for the singular solutions for an infinite wedge
geometry \cite{williams1952stress}, as Eq.\,\eqref{eq:76} is equivalent to Eq.\,\eqref{eq:29}.
More generally, for $\pi/2<\alpha<\pi$ (cases like Fig.\,\ref{fig:12}) there
is an integrable singularity with $1<t_1^*<2$, and for $0<\alpha<\pi/2$ (cases like Fig.\,\ref{fig:8})
there is no singularity because $0<t_1^*<1$.
\begin{acknowledgements}
We thank Jeremy Hoskins for a helpful discussion on numerical aspects of this work.\end{acknowledgements}

\section*{\textemdash \textemdash \textemdash \textemdash \textemdash \textemdash \textemdash{}}

\bibliographystyle{spbasic}
\addcontentsline{toc}{section}{\refname}\bibliography{Reference}

\end{document}